\documentclass[11pt]{article}

\usepackage[margin=1in]{geometry}
\usepackage[authoryear,round]{natbib}
\usepackage{graphicx}
\usepackage[hidelinks]{hyperref}

\usepackage{enumitem}
\usepackage{booktabs}
\usepackage{multirow}
\usepackage{pdflscape}
\usepackage{amssymb,amsthm,amsmath}
\usepackage{caption}
\usepackage{subcaption}
\usepackage{url}
\usepackage[english]{babel}
\usepackage[utf8]{inputenc}
\usepackage{algpseudocode,algorithm,algorithmicx}
\usepackage{makecell}
\usepackage[dvipsnames]{xcolor}
\graphicspath{{Figures/}}

\newtheorem{theorem}{Theorem}

\newtheorem{definition}[theorem]{Definition}

\newcommand{\vecU}{\mathbf{U}}
\newcommand{\vecF}{\mathbf{F}}
\newcommand{\vecx}{\mathbf{x}}

\newcommand{\vecq}{\mathbf{q}}
\newcommand{\grad}{\nabla_\mathbf{x}}
\newcommand{\R}{\mathbb{R}}

\makeatletter
\renewcommand\paragraph{\@startsection{paragraph}{4}{\z@}%
            {-2.5ex \@plus -1ex \@minus -.25ex}%
            {1.25ex \@plus .25ex}%
            {\normalfont\normalsize\itshape}}
\makeatother

\usepackage{adjustbox}

\begin{document}

\title{Efficient Weak-Entropy PINN for Solving Hyperbolic Conservation Laws}
\author{Qi Gao$^{1}$ \and Kuang Huang$^{2}$ \and Xuan Di$^{1,3}$\thanks{Corresponding author. Tel.: +1 212 853 0435; Email: \texttt{sharon.di@columbia.edu}}}
\date{\small
$^{1}$Department of Civil Engineering and Engineering Mechanics, Columbia University\\
$^{2}$Department of Mathematics, The Chinese University of Hong Kong\\
$^{3}$Center for Smart Cities, Data Science Institute, Columbia University}

\maketitle

\begin{abstract}
In recent years, neural networks have significantly advanced numerical solutions of partial differential equations (PDEs). However, solving PDEs with discontinuous solutions, such as hyperbolic conservation laws, remains challenging for neural network-based methods such as physics-informed neural networks (PINNs). Existing methods often rely on strong prior assumptions such as knowledge of discontinuity locations, or they introduce artificial smoothing terms that degrade accuracy. However, accurately solving these conservation laws and predicting the formation and propagation of discontinuities in solutions is crucial in many practical applications, including gas dynamics and traffic flow modeling. In this paper, we introduce a novel Weak-Entropy PINN (WEPINN) framework for hyperbolic conservation laws with discontinuous solutions. The method enforces the governing equations in their weak (integral) formulation and incorporates the entropy condition to select the physically admissible solution, while employing the discrete fast Fourier transform (DFFT) for efficient numerical integration. Our method is tested through extensive numerical experiments on a variety of scalar conservation laws and systems of conservation laws in one and two dimensional spaces. These experiments demonstrate that our method can accurately resolve sharp discontinuities while effectively capturing interactions between multiple shock and rarefaction waves.

\end{abstract}


\section{Introduction}

Recent years have seen an exponential growth in the use of neural networks for solving partial differential equations (PDEs), e.g., physics-informed neural networks (PINNs) \citep{raissi2019physics}. Despite their success in a variety of applications, PINNs require that solutions to the PDEs possess some level of smoothness and face significant challenges with discontinuous solutions, which often appear in hyperbolic conservation laws:
\begin{align}\label{eq:HCL}
    \partial_t \vecU(t,\vecx) + \grad \cdot \vecF(\vecU(t,\vecx))=0,
\end{align}
where $\vecU: (0,\infty) \times \R^d \to \R^m$ denotes the vector of conserved 
quantities, $t > 0$ is time, $\vecx = (x_1, \dots, x_d)^\top \in \R^d$ is the spatial 
coordinate in $d$ dimensions, and $\vecF: \R^m \to \R^{m\times d}$ is the flux function with columns 
$F_1(\vecU), \dots, F_d(\vecU)$, where each column $F_j:\R^m\to\R^m$ corresponds to the $j$-th flux component. The divergence operator acts as $\grad \cdot \vecF(\vecU) = \sum_{j=1}^d \partial_{x_j} F_j(\vecU)$.

These conservation laws are widely used in various scientific and engineering disciplines, such as (inviscid) Burgers' equation and Euler equations for gas dynamics \citep{dafermos1983hyperbolic}, or Lighthill-Whitham-Richards (LWR) model for traffic flow \citep{lighthill1955kinematic,MR75522}.
Solutions to these equations often develop discontinuities even from smooth initial conditions, giving rise to critical phenomena such as shock waves in gas propagation or traffic congestion. 
Traffic congestion can be understood as a discontinuity in vehicle density when vehicles upstream brake abruptly, while entering a jammed region downstream.

The standard PINN framework \citep{raissi2019physics} embeds governing equations in their strong form directly into the loss function by minimizing pointwise PDE 
residuals, e.g.,
\begin{equation}\label{eq:diffpinn_pde_loss}
\mathcal{L}_{\text{PDE}} = \frac{1}{N} \sum_{i=1}^{N} \left\| \partial_t \hat{\vecU}(t_i, \vecx_i) + \grad \cdot \vecF(\hat{\vecU}(t_i, \vecx_i)) \right\|^2,
\end{equation}
computed via automatic differentiation, where $\hat{\vecU}$ is the neural network approximation and $\{(t_i,\vecx_i)\}_{i=1}^N$ is a set of collocation points. Because the residual is evaluated pointwise through differentiation, we refer to this method as \textbf{Diff-PINN}. This approach inherently relies on solution smoothness and therefore struggles with discontinuous solutions that arise in hyperbolic conservation laws, typically producing overly smooth profiles or spurious oscillations near discontinuities \citep{Huang2023OnTL}.
 
When solving hyperbolic conservation laws with discontinuous solutions, one workaround to enable the use of pointwise residuals in the Diff-PINN
is to incorporate artificial viscosity for regularization, i.e., to consider
\[ \partial_t \vecU(t,\vecx) +\grad \cdot \vecF(\vecU(t,\vecx))=\varepsilon \Delta_{\vecx} \vecU(t,\vecx), \]
for a suitable regularization parameter $\varepsilon>0$, where $\Delta_{\vecx}=\sum_{j=1}^d \partial_{x_j}^2$ is the Laplacian operator. Although this makes solutions smooth, it also permanently smooths out the discontinuities, producing smeared shock profiles whose width is controlled by $\varepsilon$, rather than reproducing the sharp discontinuities present in the true solution.
In contrast, classical numerical methods such as the Lax-Friedrichs and Godunov schemes also introduce numerical viscosity, but this viscosity vanishes as the mesh is refined, and 
convergence to the physically correct solution is guaranteed by the mathematical framework of weak formulations and entropy conditions \citep{leveque2002finite}. This well-established  framework has inspired several weak-form PINN approaches.

Weak PINN (\textbf{WPINN}) \citep{de2024wpinns} incorporates residuals computed from the weak formulation and the entropy condition into the loss function. The key feature of WPINN is that it adopts an adversarial framework in which two neural networks are trained simultaneously: one to approximate the solution and the other to generate test functions. These networks are trained with opposite objectives: the solution network minimizes the residual loss, while the test function network maximizes it. This adversarial setup is designed to identify the most informative test functions without exhaustively searching the test function space, thereby improving efficiency. \cite{chaumet2022improving} proposed several modifications to the WPINN framework, including dual norm computation for training efficiency, as well as extensions to general boundary conditions and systems of conservation laws. \cite{zhou2025weak} further extended WPINN to complex geometries. However, as noted in \cite{de2024wpinns, chaumet2022improving}, the adversarial training of WPINN faces stability challenges in practice. The the test function network may require frequent re-initialization during training, and the predictions can be sensitive to hyperparameter choices.

Variational PINN (\textbf{VPINN}) \citep{kharazmi2019variational, kharazmi2021hp} is another variant of PINNs that leverages integration by parts to transfer derivatives from the predicted solution to test functions. The VPINN formulates its loss function as the sum of residuals from the resulting integral equations. Subsequent work analyzed the role of quadrature rules and test function selection in VPINNs \citep{Berrone2021VariationalPI}, and proposed a robust variant based on discrete dual norm minimization \citep{Rojas2023RobustVP}. Although it is not specifically developed for hyperbolic conservation laws, one of its residual terms naturally corresponds to a weak formulation. However, VPINN does not enforce the entropy condition, which is essential for selecting physically admissible solutions among multiple weak solutions.

Beyond weak-form approaches, a variety of other strategies have been proposed to address discontinuities when using PINNs.
For instance, \cite{jagtap2020extended, jagtap2020conservative, lorin2024non} adopt domain decomposition by dividing the solution domain into sub-regions separated at the shock locations. However, this approach requires prior knowledge of shock locations, limiting its applicability.
More complex training pipelines \citep{wang2025solving} have also been explored, in which a PINN is iteratively trained to identify shock locations, followed by domain decomposition into subdomains that contain no discontinuities. Independent PINNs are then trained within each subdomain, but this approach becomes increasingly impractical when multiple discontinuities interact or merge over time, as the subdomain configuration must be continuously updated.
Gradient-weighting strategies \citep{liu2024discontinuity, liang2024continuous, ghoreishi2026physics} suppress large residual gradients near shocks to mitigate training instability. While effective for problems with stationary or isolated discontinuities, their performance degrades when shocks are moving or multiple discontinuities interact. Network architecture improvements \citep{lei2025KAN, lei2025discontinuity} provide another approach by incorporating Kolmogorov-Arnold network (KAN) layers with explicit high-frequency feature embeddings to better capture shock waves. However, these methods still require additional regularization such as adaptive artificial viscosity.

Alternative reformulations include a lift-and-embed method \citep{sun2024lift} that lifts the problem to a higher dimension where solutions become smooth, though it relies on known initial locations of discontinuities. In \cite{zhang2022implicit}, an implicit form is used to design the loss function, thereby avoiding differentiation of the predicted solution. Subsequent work \citep{zeng2025clinn} extends this framework with solution bounds and Rankine-Hugoniot jump conditions to improve shock detection and the accuracy of predicted shock speeds. However, the implicit form does not enforce the entropy condition that is necessary for selecting physically admissible solutions. In \cite{su2025spike}, a Stable Physics-Informed Kernel Evolution (SPIKE) method is proposed. This method employs reproducing kernel representations of solutions, which are equivalent to two-layer neural networks, together with Tikhonov regularization. However, its applicability to multidimensional problems remains to be explored.
\cite{cai2022least, cai2023least} proposed a least-squares ReLU neural network (LSNN) that leverages finite-volume discretization of conservation laws, though its extension to systems of conservation laws has not been demonstrated.

The aforementioned methods either require prior knowledge of discontinuity locations, lack essential constraints such as the entropy condition, introduce artificial smoothing that degrades accuracy near discontinuities, or suffer from training instability and inefficiency.
These challenges become particularly severe in regimes involving complex nonlinear wave dynamics, including shock formation from initially smooth profiles, shock merging, and shock-rarefaction interactions.

To overcome these challenges, we propose a novel PINN framework for hyperbolic conservation laws that rigorously integrates weak formulations and entropy conditions while ensuring stable and efficient training.
This mathematically grounded foundation mirrors the classical theory and numerical methods for hyperbolic conservation laws, making our approach applicable to general systems of conservation laws in multiple space dimensions. We refer to this method as the Weak-Entropy PINN (\textbf{WEPINN}), which addresses the challenges of discontinuous solutions through three key components:
\begin{itemize}
\item the weak formulation of the conservation laws as integral equations, which transfers all spatial and temporal derivatives onto smooth test functions (e.g., trigonometric polynomials) and inherently accommodates discontinuities without artificial smoothing or prior knowledge of their locations;
\item enforcement of the entropy condition as integral inequalities to select the physically admissible solution, enabling the method to correctly determine whether an initial discontinuity evolves into a shock or a rarefaction wave;
\item use of the discrete fast Fourier transform (DFFT) for efficient and accurate evaluation of the integrals appearing in both the weak formulation and the entropy condition, significantly reducing computational cost compared to standard numerical integration.
\end{itemize}

A key feature of the proposed WEPINN is that it pre-selects a set of test functions and minimizes the residuals over all of them simultaneously via least-squares optimization, avoiding the adversarial training employed in WPINN and leading to stable and efficient training. 
We demonstrate the effectiveness of WEPINN through extensive numerical experiments covering one-dimensional scalar conservation laws (the linear advection equation, Burgers' equation, and the LWR traffic model), one-dimensional system of conservation laws (the compressible Euler 
equations), and two-dimensional scalar conservation laws (the 2D Burgers' equation), under a variety of initial and boundary conditions. Our experiments systematically evaluate a wide range of nonlinear wave phenomena, including shock formation from smooth initial conditions, propagation of shock and rarefaction waves, shock merging, and shock-rarefaction interactions.
Furthermore, we introduce two shock-aware evaluation metrics, the shock detection rate and shock position accuracy, which complement traditional $\mathrm{L}^p$ error measures by directly assessing the presence and location of discontinuities predicted by a model.

The remainder of this paper is organized as follows. Section~\ref{sec:methods} introduces the mathematical framework of weak entropy solutions and presents the detailed design of the proposed Weak-Entropy PINN (WEPINN). 
Section~\ref{sec:results} provides extensive numerical experiments and performance evaluations. 
Section~\ref{sec:ablation} presents ablation studies that examine the role of the entropy condition and the choice of test functions in WEPINN. Section~\ref{sec:conclusion} concludes the paper and discusses potential future directions.

\section{Methods}\label{sec:methods}

In this section, we present the Weak-Entropy PINN (WEPINN) framework. We begin by reviewing the theory of weak entropy solutions to hyperbolic conservation laws, which provides the mathematical foundation of our method, and then describe the design of the WEPINN.

\subsection{Weak entropy solutions to hyperbolic conservation laws}
We consider the hyperbolic conservation law \eqref{eq:HCL} on a hypercube $\Omega=[a,b]^d\subset\R^d$ for some $a<b$, equipped with the initial condition
\begin{equation}\label{eq:ini_cond}
    \vecU(0,\vecx)=\vecU_0(\vecx) \quad \text{for} \quad \vecx\in\Omega,
\end{equation}
and periodic boundary conditions. In this work, we address non-periodic boundary conditions only for one-dimensional conservation laws (cf.~Section~\ref{sec:boundary_cond}).

For discontinuous solutions to \eqref{eq:HCL}, the derivatives $\partial_t\vecU$ and $\grad\cdot\vecF(\vecU)$ in \eqref{eq:HCL} do not exist at the discontinuities, so that the equation fails to describe the propagation of the discontinuities. This necessitates considering weak solutions defined by the following integral equation:
\begin{align}
    [\text{Weak}] \quad &\int_0^\infty\int_\Omega \vecU(t,\vecx)\partial_t\varphi(t,\vecx) + \vecF(\vecU(t,\vecx))\cdot\grad\varphi(t,\vecx) \,d\vecx dt + \int_\Omega \vecU_0(\vecx) \varphi(0,\vecx)\,d\vecx = 0,
    \label{eq:weak}
\end{align}
where $\varphi:[0,\infty)\times\Omega\to\R$ is any spatially periodic, temporally compactly supported, and continuously differentiable test function.

The integral equation \eqref{eq:weak}, known as the \emph{weak formulation} of equation \eqref{eq:HCL}, is obtained by multiplying \eqref{eq:HCL} by the test function $\varphi$ and integrating by parts. The weak formulation \eqref{eq:weak} imposes no smoothness condition on $\vecU$ and yields the Rankine-Hugoniot jump condition \citep{dafermos1983hyperbolic} that capture the propagation of discontinuities. However, weak solutions defined by \eqref{eq:weak} are generally non-unique; in particular, multiple weak solutions can correspond to the same initial data, differing in how discontinuities form. To select physically admissible solutions, one crucial criterion is the entropy condition:
\begin{align}
    [\text{Entropy}] \quad &\int_0^\infty\int_\Omega \eta(\vecU(t,\vecx))\partial_t\varphi(t,\vecx) + \vecq(\vecU(t,\vecx)) \cdot \grad\varphi(t,\vecx) \,d\vecx dt + \int_\Omega \eta(\vecU_0(\vecx)) \varphi(0,\vecx)\,d\vecx \geq0, 
    \label{eq:entropy}
\end{align}
where $\eta:\R^m\to\R$ and $\vecq:\R^m\to\R^d$ form an \emph{entropy-entropy flux pair} satisfying 
\begin{equation}\label{eq:entropy_pair}
    \frac{\partial q_j}{\partial u_k} = \sum_{l=1}^m \frac{\partial\eta}{\partial u_l}\frac{\partial F_{j,l}}{\partial u_k} \quad \text{for} \quad j=1,\cdots,d, \quad k=1,\cdots,m,
\end{equation}
where $u_k$ denotes the $k$-th component of $\vecU$, $q_j$ the $j$-th component of $\vecq$, and $F_{j,l}$ the $l$-th component of the $j$-th flux $F_j$ of $\vecF$.
When the solution $\vecU$ is smooth, the chain rule and \eqref{eq:entropy_pair} yield a conservation law for the entropy:
\begin{align*}
    \partial_t \eta(\vecU) + \grad\cdot \vecq(\vecU) = 0.
\end{align*}
When discontinuities form, however, this equation no longer holds.
The entropy condition \eqref{eq:entropy} replaces it with integral inequalities over all nonnegative test functions $\varphi$, characterizing the dissipation of the entropy $\eta(\vecU)$ and the irreversibility of the time evolution of $\vecU$.

The weak formulation \eqref{eq:weak} and the entropy condition \eqref{eq:entropy} together give the following definition of \emph{weak entropy solutions}.
\begin{definition}\label{def:weak_entropy_sol}
Suppose that $\mathcal{E}$ is a family of entropy-entropy flux pairs $(\eta,\vecq)$ that satisfy condition \eqref{eq:entropy_pair}, and let $\mathcal{T}$ denote the space of functions $\varphi\in\mathrm{C}^1([0,\infty) \times \Omega)$ that are spatially periodic on $\Omega$ and compactly supported in time. A function $\vecU\in\mathrm{L}^\infty((0,\infty)\times\Omega)\cap\mathrm{C}((0,\infty);\mathrm{L}^1(\Omega))$ is called a weak entropy solution to the hyperbolic conservation law \eqref{eq:HCL} with the initial condition \eqref{eq:ini_cond} if it satisfies the weak formulation \eqref{eq:weak} for all test functions $\varphi\in\mathcal{T}$ and the entropy condition \eqref{eq:entropy} for all nonnegative test functions $\varphi \in \mathcal{T}$ and every $(\eta, \vecq) \in \mathcal{E}$.
\end{definition}

For scalar conservation laws ($m=1$), any $\mathrm{C}^1$ convex function $\eta:\R\to\R$ serves as an entropy, with the corresponding entropy flux $\vecq:\R\to\R^d$ defined component-wise by  $q_j'(\xi) = \eta'(\xi) F_j'(\xi)$ for $j = 1, \dots, d$. The existence and uniqueness of weak entropy solutions in this case are well established for any spatial dimension $d$, following the classical theory of Kru\v{z}kov \citep{kruvzkov1970first}.

\begin{theorem}\label{thm:well-posedness}
Suppose that $m=1$. For any initial data $\vecU_0\in\mathrm{L}^\infty(\Omega)$, there exists a unique weak entropy solution $\vecU\in\mathrm{L}^\infty((0,\infty)\times\Omega)\cap\mathrm{C}((0,\infty);\mathrm{L}^1(\Omega))$ to the scalar conservation law \eqref{eq:HCL} in the sense of Definition~\ref{def:weak_entropy_sol}, where the family $\mathcal{E}$ consists of all entropy-entropy flux pairs $(\eta, \vecq)$ with $\mathrm{C}^1$ convex $\eta$.
\end{theorem}

For systems of conservation laws ($m>1$), no generic construction of entropy-entropy flux pairs is available. However, certain systems arising from gas and fluid dynamics admit natural entropy structures. For example, the compressible Euler equations possess entropy pairs derived 
from thermodynamic entropy, as detailed in Section~\ref{sec:vector-1d}. The existence and uniqueness of weak entropy solutions for general systems remain open problems 
\citep{dafermos1983hyperbolic,Dafermos2023} and are beyond the scope of this work. 

\begin{figure}[htbp]
    \centering
    \includegraphics[width=0.8\linewidth]{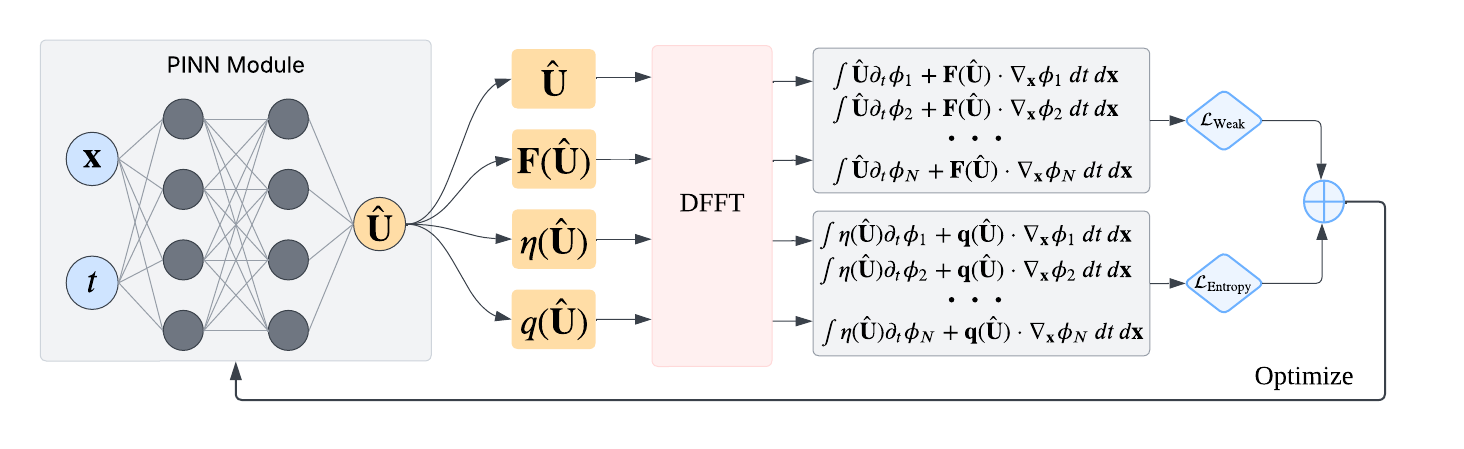}
    \caption{Pipeline of the proposed WEPINN}
    \label{fig:pipeline}
\end{figure}

\subsection{Neural network approximation to solutions}

We approximate solutions to the hyperbolic conservation law \eqref{eq:HCL} using a neural network
\begin{equation}
(t,\mathbf{x})\in\R^{d+1} \longmapsto \hat{\mathbf{U}}(t,\mathbf{x};\theta)\in\R^m,
\end{equation}
where $\theta$ denotes the collection of all trainable parameters (weights and biases) of the network. The solutions may be discontinuous but remain of bounded variation, a regularity class that neural networks can approximate to arbitrary accuracy \citep{DeRyck2024NumericalAO}.

In the Diff-PINN framework, the neural network is trained by minimizing a loss function that enforces the governing equation \eqref{eq:HCL} in its strong form and the initial and boundary conditions. The loss function is defined as
\begin{equation*}
    \mathcal{L}_{\text{Strong}} = \mathcal{L}_{\text{PDE}} + \mathcal{L}_{\text{IC}} + \mathcal{L}_{\text{BC}},
\end{equation*}
where $\mathcal{L}_{\text{PDE}}$ denotes the pointwise PDE residual as in \eqref{eq:diffpinn_pde_loss}, and $\mathcal{L}_{\text{IC}}$ and $\mathcal{L}_{\text{BC}}$ penalize deviations from the prescribed initial and boundary conditions at collocation points sampled on the initial time slice and along the domain boundary, respectively.

For hyperbolic conservation laws, the true solution $\vecU$ may develop discontinuities, making $\mathcal{L}_{\text{Strong}}$ an inappropriate training objective.
In our WEPINN framework, we replace $\mathcal{L}_{\text{Strong}}$ with 
a weak-entropy loss based on \eqref{eq:weak}--\eqref{eq:entropy}, as described in the following subsection.

\subsection{Weak-Entropy PINN}\label{eq:method_wepinn}

The proposed WEPINN is built on the notion of weak entropy solutions defined in Definition~\ref{def:weak_entropy_sol}. As illustrated in Figure~\ref{fig:pipeline}, we employ a standard PINN architecture to represent the approximate solution as a neural network: 
$(t, \vecx) \mapsto \hat{\vecU}(t, \vecx; \theta)$, and design a weak-entropy loss by discretizing \eqref{eq:weak}--\eqref{eq:entropy} to train the neural network.

We discretize~\eqref{eq:weak}–\eqref{eq:entropy} on a finite time horizon $[0,T]$ where $T>0$ is the terminal time of interest, using a spatial mesh $\{\vecx_j\}_{j\in\Omega_d}$ on $\Omega$ and a temporal mesh $\{t_k\}_{k=0}^{N_t}$ on $[0,T]$, with $\{\beta_j\}_{j\in\Omega_d}$ and $\{\alpha_k\}_{k=0}^{N_t}$ the corresponding quadrature weights.
That is, for any function $h:[0,T]\times\Omega\to\R$, the integral $\int_0^T\int_\Omega h(t,\vecx)\, d\vecx dt$ is discretized as
\[ \sum_{j\in\Omega_d}\sum_{k=0}^{N_t} \alpha_k\beta_j h(t_k,x_j). \]
In our implementation, both $\{\alpha_k\}_{k=0}^{N_t}$ and $\{\beta_j\}_{j\in\Omega_d}$ are taken as trapezoidal-rule quadrature weights on the uniform grid.

For any test function $\varphi$, we discretize the weak formulation \eqref{eq:weak} as
\[ \mathcal{I}_{t,\vecx}^\varphi + \mathcal{I}_{0,\vecx}^\varphi = 0, \]
through numerical quadratures:
\begin{align}
    \mathcal{I}_{t,\vecx}^\varphi&=\sum_{j\in\Omega_d}\sum_{k=0}^{N_t} \alpha_k\beta_j \big(\hat{\vecU}(t_k,\vecx_j;\theta) \partial_t\varphi(t_k,\vecx_j) + \vecF(\hat{\vecU}(t_k,\vecx_j;\theta)) \cdot \grad \varphi(t_k,\vecx_j) \big), \label{eq:Ixt} \\
    \mathcal{I}_{0,\vecx}^\varphi&= \sum_{j\in\Omega_d} \beta_j \vecU_0(\vecx_j)\, \varphi(0,\vecx_j), \label{eq:Ix}
\end{align}

Then, for any nonnegative test function $\varphi$, we discretize the entropy condition \eqref{eq:entropy} as
\[ \mathcal{J}_{t,x}^\varphi+\mathcal{J}_{0,x}^\varphi\geq 0,\]
through numerical quadratures:
\begin{align}
    \mathcal{J}_{t,x}^\varphi=&\sum_{j\in\Omega_d}\sum_{k=0}^{N_t} \alpha_k\beta_j \big(\eta(\hat{\vecU}(t_k,\vecx_j;\theta)) \partial_t \varphi(t_k,\vecx_j) + \vecq(\hat{\vecU}(t_k,\vecx_j;\theta)) \cdot \grad\varphi(t_k,\vecx_j) \big), \label{eq:Jxt} \\
    \mathcal{J}_{0,\vecx}^\varphi=& \sum_{j\in\Omega_d} \beta_j \eta(\vecU_0(\vecx_j))\varphi(0,\vecx_j). \label{eq:Jx}
\end{align}

With the numerical discretization of the weak formulation and the entropy condition above,
our weak-entropy loss is defined as
\begin{align*}
\mathcal{L}_{\mathrm{Weak-Entropy}} = \mathcal{L}_{\mathrm{Weak}}+\mathcal{L}_{\mathrm{Entropy}} + \mathcal{L}_{\mathrm{IC}},
\end{align*}
where
\begin{align*}
\mathcal{L}_{\mathrm{Weak}} &= \frac1N\sum_{n=1}^N \Big\|\mathcal{I}_{t,\vecx}^{\varphi_n}+\mathcal{I}_{0,\vecx}^{\varphi_n}\Big\|, \\
\mathcal{L}_{\mathrm{Entropy}} &= \frac1N\sum_{n=1}^N \text{ReLU}\left(-\mathcal{J}_{t,\vecx}^{\varphi_n}-\mathcal{J}_{0,\vecx}^{\varphi_n},\,0\right), \\
\mathcal{L}_{\text{IC}} &= \frac{1}{N_0} \sum_{i=1}^{N_0} 
\left\| \hat{\vecU}(0,\vecx_i^0;\theta) - \mathbf{U}_0(\vecx_i^0) \right\|^2,
\end{align*}
with $\mathcal{I}_{t,\vecx}^{\varphi_n},\mathcal{I}_{0,\vecx}^{\varphi_n},\mathcal{J}_{t,\vecx}^{\varphi_n},\mathcal{J}_{0,\vecx}^{\varphi_n}$ as defined in equations \eqref{eq:Ixt}--\eqref{eq:Jx}, $\{\varphi_n\}_{n=1}^N$ is a set of pre-selected test functions, and $\{\vecx_i^0\}_{i=1}^{N_0}$ a set of collocation points sampled on the initial time slice $t=0$.

\paragraph{Test functions.}
The weak formulation in principle requires evaluating the residuals in \eqref{eq:weak} against all test functions $\varphi\in\mathcal{T}$ as in Definition~\ref{def:weak_entropy_sol}. Since this is computationally infeasible, we approximate $\mathcal{T}$ by a finite-dimensional subspace spanned by a representative set of test functions and evaluate the residuals only on this subset.

We adopt a trigonometric basis on a uniform grid, with a tensor-product structure in space and time. The tensor-product structure allows the high-dimensional summations in the weak-entropy loss to factorize along each coordinate direction, and combined with a uniform grid and trigonometric basis functions, enables an efficient implementation through the \emph{discrete fast Fourier transform} (DFFT), as described below. The choice of the trigonometric basis over alternatives such as orthogonal polynomials is supported by an ablation study in Section~\ref{sec:test-function-abl}, with additional discussion in \ref{sec:test-function-explan}.

Concretely, each test function in our finite set $\{\varphi_n\}_{n=1}^N$ is indexed by a multi-index
\[
\nu(n) = (\nu_0(n), \nu_1(n), \dots, \nu_d(n)),
\]
where $\nu_0(n) \in \{0,\dots,N_q\}$ selects the temporal basis and $\nu_r(n) \in \{0,\dots,N_p\}$ selects the spatial basis in the $r$-th spatial dimension. The total number of test functions is $N = (N_q + 1)(N_p + 1)^d$, and each test function factorizes as
\begin{equation}
\varphi_n(t,\mathbf{x})
=
q_{\nu_0(n)}(t)\prod_{r=1}^{d} p_{r,\nu_r(n)}(x_r).
\end{equation}
We use two sets of basis functions: one for the weak formulation, where test functions only need to be smooth, and one for the entropy condition, where test functions must additionally satisfy $\varphi\geq0$.

Suppose that the spatial domain $\Omega=[0,1]^d$ and the time horizon is $[0,T]$. For the weak formulation, we adopt standard trigonometric bases. The spatial basis per dimension consists of sine-cosine pairs with frequencies $m=1,\dots,N_p/2$, where $N_p$ is assumed even.
For each spatial dimension $r=1,\dots,d$:
\begin{equation}
p_{r,0}(x_r)=1, \quad 
p_{r, 2m-1}(x_r) = \sin(2\pi m x_r), \quad
p_{r, 2m}(x_r) = \cos(2\pi m x_r), \quad m=1,\dots,N_p/2.
\end{equation}
For the temporal basis, we introduce a window function $w(t)=T-t$ to ensure $\varphi(T,\vecx)=0$.
The temporal basis uses sine-cosine pairs with frequencies $m/2$ for $m=1,\dots,N_q/2$, where $N_q$ is assumed even:
\begin{equation}
q_0(t)=w(t), \quad
q_{2m-1}(t) = w(t)\sin(\pi m t), \quad
q_{2m}(t) = w(t)\cos(\pi m t), \quad m=1,\dots,N_q/2.
\end{equation}

For the entropy condition, we modify the above bases to guarantee non-negativity. The spatial basis uses four nonnegative combinations per frequency $m=1,\dots,N_p/4$, where $N_p$ is assumed divisible by 4:
\begin{align}
&p_{r, 0}(x_r) = 1, \quad
p_{r, 4m-3}(x_r) = 1 + \sin(2\pi m x_r), \quad
p_{r, 4m-2}(x_r) = 1 + \cos(2\pi m x_r), \\
&p_{r, 4m-1}(x_r) = 1 - \sin(2\pi m x_r), \quad
p_{r, 4m}(x_r) = 1 - \cos(2\pi m x_r), \quad m=1,\dots,N_p/4,
\end{align}
and the temporal basis uses four nonnegative combinations per frequency $m/2$ for $m=1,\dots,N_q/4$, where $N_q$ is assumed divisible by 4, with the window function $w(t)$:
\begin{align}
&q_{0}(t) = w(t), \quad
q_{4m-3}(t) = w(t)\big(1 + \sin(\pi m t)\big), \quad
q_{4m-2}(t) = w(t)\big(1 + \cos(\pi m t)\big), \\
&q_{4m-1}(t) = w(t)\big(1 - \sin(\pi m t)\big), \quad
q_{4m}(t) = w(t)\big(1 - \cos(\pi m t)\big), \quad m=1,\dots,N_q/4.
\end{align}

\paragraph{DFFT implementation.}
The dominant computational cost of evaluating the weak-entropy loss lies in the residuals $\mathcal{I}_{t,\mathbf{x}}^{\varphi_n}$ and $\mathcal{J}_{t,\mathbf{x}}^{\varphi_n}$, which involve summations over the full space-time grid. For high-frequency test functions, direct evaluation of these summations becomes both time- and memory-intensive.

Our choice of trigonometric test functions on a uniform grid allows these summations to be evaluated via the DFFT. The resulting computational complexity is reduced from $O(N_t^2 N_x^{2d})$ for direct numerical quadrature to $O(N_t N_x^d \log N_t \log^d N_x)$, where $N_t$ is the number of temporal grid points and $N_x$ is the number of grid points per spatial dimension. Implementation details and empirical validation of the speedup are provided in \ref{sec:efficiency}.

\subsection{Boundary conditions}\label{sec:boundary_cond}

The weak‑entropy loss defined in Section~\ref{eq:method_wepinn} assumes periodic boundary conditions. With the chosen periodic test functions, the boundary terms in the weak formulation \eqref{eq:weak} and the entropy condition \eqref{eq:entropy} vanish identically, and no additional boundary loss is required in WEPINN.

For one-dimensional problems, the boundary of the domain $\Omega=[a,b]$ consists of two points $x=a$ and $x=b$. In this case, we also consider the Dirichlet boundary condition:
\[ \vecU(t,a) = \vecU_0(a) \quad \text{and} \quad \vecU(t,b)=\vecU_0(b) \quad \text{for } t\in[0,T]. \]
This boundary condition is interpreted as extending the problem to the real line $\R$ by the following constant extension of solutions:
\[ \tilde{\vecU}_0(x) = \begin{cases}
    \vecU_0(a), \quad &x\leq a; \\
    \vecU_0(x), \quad &a< x <b; \\
    \vecU_0(b), \quad &x\geq b;
\end{cases} \qquad \tilde{\vecU}(t,x) = \begin{cases}
    \vecU_0(a), \quad &x\leq a; \\
    \vecU(t,x), \quad &a< x <b; \\
    \vecU_0(b), \quad &x\geq b;
\end{cases}\]
Then, we can study the initial value problem \eqref{eq:HCL} for $\tilde{\vecU}$ with the initial condition $\tilde{\vecU}(0,x)=\tilde{\vecU}_0(x)$ for $x\in\R$.
Provided that all waves propagate within the domain $[a,b]$ over the time horizon $[0,T]$, the extended solution coincides with the original one on $[a,b]$, i.e., $\vecU(t,x)=\tilde{\vecU}(t,x)$ for all $t\in[0,T]$ and $x\in[a,b]$.
This allows us to consider Riemann problems in one dimension.

For the Dirichlet boundary condition, the weak formulation and the entropy condition become
\begin{align*}
    [\text{Weak}] \quad &\int_0^T\int^b_a \vecU(t,x)\partial_t\varphi(t,x) + \vecF(\vecU(t,x))\partial_x\varphi(t,x) \,dx dt + \int^b_a \vecU_0(x) \varphi(0,x)\,dx \\
   & + \int_0^T \vecF(\vecU_0(a))\varphi(t, a)\,dt - \int_0^T \vecF(\vecU_0(b))\varphi(t, b) \,dt = 0, \\
    [\text{Entropy}] \quad &\int_0^T\int^b_a \eta(\vecU(t,x))\partial_t\varphi(t,x) + \vecq(\vecU(t,x))\partial_x\varphi(t,x) \,dx dt + \int^b_a \eta(\vecU_0(x)) \varphi(0,x)\,dx \\
    & + \int_0^T \vecF(\vecU_0(a))\varphi(t, a)\,dt - \int_0^T \vecF(\vecU_0(b))\varphi(t, b)\,dt \geq0, 
\end{align*}
respectively. Correspondingly, the weak and entropy losses are updated as follows:
\begin{align*}
\mathcal{L}_{\mathrm{Weak}} &= \frac1N\sum_{n=1}^N \Big\|\mathcal{I}_{t,x}^{\varphi_n}+\mathcal{I}_{0,x}^{\varphi_n}+\mathcal{I}_{\mathrm{BC}}^{\varphi_n}\Big\|, \\
\mathcal{L}_{\mathrm{Entropy}} &= \frac1N\sum_{n=1}^N \text{ReLU}\left(-\mathcal{J}_{t,x}^{\varphi_n}-\mathcal{J}_{0,x}^{\varphi_n}-\mathcal{J}_{\mathrm{BC}}^{\varphi_n},\ 0\right),
\end{align*}
where $\mathcal{I}_{t,x}^{\varphi_n},\mathcal{I}_{0,x}^{\varphi_n},\mathcal{J}_{t,x}^{\varphi_n},\mathcal{J}_{0,x}^{\varphi_n}$ are as in Section~\ref{eq:method_wepinn}, and the boundary terms $\mathcal{I}_{\mathrm{BC}}^{\varphi_n}$ and $\mathcal{J}_{\mathrm{BC}}^{\varphi_n}$ are defined as:
\begin{align*}
    \mathcal{I}_{\mathrm{BC}}^{\varphi_n}&= \sum_{k=0}^{N_t} \alpha_k\vecF(\vecU_0(a)) {\varphi_n}(t_k,a) - \sum_{k=0}^{N_t} \alpha_k\vecF(\vecU_0(b)) {\varphi_n}(t_k,b),\\
    \mathcal{J}_{\mathrm{BC}}^{\varphi_n}&= \sum_{k=0}^{N_t} \alpha_k \vecq(\vecU_0(a)) {\varphi_n}(t_k,a) - \sum_{k=0}^{N_t} \alpha_k \vecq(\vecU_0(b)) {\varphi_n}(t_k,b).
\end{align*}

\section{Numerical Experiments and Results}\label{sec:results}
In this section, we evaluate the performance of the proposed method through comprehensive numerical experiments. Specifically, we consider scalar conservation laws and systems of conservation laws in one and two spatial dimensions, as summarized in Table~\ref{tab:pde-benchmarks}.

\begin{table}[H]
    \centering
    \renewcommand{\arraystretch}{1.15}
    \begin{tabular}{|c|c|c|}
        \hline
        \textbf{Problem class} & \textbf{One-dimension} & \textbf{Two-dimension} \\ \hline
        Scalar conservation laws 
        & 
        \makecell[l]{Linear advection equation \\ Burgers' equation \\ LWR traffic model \\ \emph{(Section~\ref{sec:scalar-1d})}}
        &
        \makecell[l]{Burgers' equation \\ \emph{(Section~\ref{sec:scalar-2d})}} 
        \\ \hline
        Systems of conservation laws 
        &
        \makecell[l]{Compressible Euler equations \\ \emph{(Section~\ref{sec:vector-1d})}}
        &
        -
        \\ \hline
    \end{tabular}
    \caption{Benchmark problems and corresponding subsections, covering both scalar conservation laws and systems of conservation laws in one and two spatial dimensions.}
    \label{tab:pde-benchmarks}
\end{table}

\subsection{Experimental setup}\label{sec:setup}

We compare the proposed WEPINN against three representative baselines that require only the governing PDEs, without prior knowledge of discontinuity locations: Diff-PINN, VPINN, and WPINN. Auxiliary techniques such as adaptive collocation sampling and domain decomposition are orthogonal to this comparison and could be combined with any of the formulations considered; therefore, we do not include them as separate baselines.

For Diff-PINN and VPINN, we adopt the same network architecture as in WEPINN, except that ReLU activations are replaced with Tanh for better performance, and we use the same number of uniformly sampled collocation points as in WEPINN. VPINN as originally proposed \citep{kharazmi2019variational, kharazmi2021hp} is not specifically designed for hyperbolic conservation laws. We adapt it by combining our weak-formulation loss with the pointwise PDE residual loss used in Diff-PINN, which corresponds to the core idea of VPINN applied to the present setting. All three methods use Adam for optimization. Full network architecture and training details are provided in Appendix C.

For WPINN, we adopt the implementation of \cite{chaumet2022improving} rather than the original version of \cite{de2024wpinns}, as it provides a TensorFlow implementation more compatible with our training environment. 
We import the initial and boundary conditions from our experimental setup into their code and evaluate the resulting solutions on a uniform mesh for downstream metric computation. The released code supports only Dirichlet boundary conditions and does not include an implementation for the compressible Euler equations. Accordingly, we report WPINN results only for the Dirichlet-boundary scalar cases, and for the Euler case we use the numbers reported in \cite{chaumet2022improving}.

To assess model performance, we evaluate both global solution accuracy and shock (discontinuity) resolution. Global accuracy is measured by the relative $\mathrm{L}^2$ error against a high-resolution reference solution computed with the Lax-Friedrichs scheme. To assess shock resolution, we first extract the shock locations from the predicted solution at each time step using a gradient-based peak detector, and then compare them against the ground-truth shock locations. This comparison yields two shock-aware metrics:
\begin{itemize}
 \item the shock detection rate (S-Rate), defined as the fraction of ground-truth shocks correctly identified by the model;
 \item the shock position accuracy (S-Acc), defined as the mean positional error between the matched predicted and ground-truth shocks. 
\end{itemize}
Together, the two metrics quantify whether a model captures the presence of shocks and how accurately it predicts their positions over a time horizon. The detailed procedure for computing these metrics is given in Algorithm~\ref{alg:metric} (Appendix~A).

\subsection{One-dimensional scalar conservation laws}\label{sec:scalar_1d}
\label{sec:scalar-1d}
In this subsection, we consider the one-dimensional scalar conservation law
\begin{equation} \label{eq:scalar-1d}
    \partial_t u(t,x) + \partial_x f(u(t,x)) = 0,
    \quad (t,x) \in [0,1] \times [0,1],
\end{equation}
subject to the initial condition \(u(0,x) = u_0(x)\) for $x\in[0,1]$ and either the Dirichlet boundary condition described in Section~\ref{sec:boundary_cond} or periodic boundary conditions. 

We consider three representative examples with different flux functions $f$:
\begin{itemize}
    \item \textbf{Linear advection equation}, with \(f(u) = cu\) for a constant advection speed \(c\).
    \item \textbf{(inviscid) Burgers' equation}, with \(f(u) = \frac12u^2\), a classical benchmark in gas dynamics.
    \item \textbf{LWR traffic model}, with \(f(u) = u(1-u)\), where \(u\) denotes traffic density.
\end{itemize}

For each combination of governing equation and boundary condition, we randomly generate 30 initial conditions drawn from several initial condition classes designed to cover a range of solution behaviors. Under the Dirichlet boundary condition, the classes are \emph{Sigmoid} (smooth transitions between two states) and \emph{Riemann} (piecewise-constant two- or three-state profiles). Under the periodic boundary condition, the classes are \emph{Fourier} (smooth trigonometric series), \emph{Trig} (high-frequency sinusoidal wave), \emph{Bell} (smooth single-peaked profiles), and \emph{PWC} (piecewise-constant profiles). Full generation procedures are given in \ref{sec:icbc}. Each model is trained independently on each initial condition, and the reported metrics are averaged across all trials.

\begin{table}[htbp]
\centering
\scriptsize
\setlength{\tabcolsep}{4pt}
\renewcommand{\arraystretch}{1.2}
\begin{tabular}{|c|c|ccc|ccc|}
\hline
\scriptsize{\textbf{Problem}} & \textbf{Method} & \multicolumn{3}{c|}{\textbf{Riemann}} & \multicolumn{3}{c|}{\textbf{Sigmoid}} \\
\cline{3-8}
&& $\mathrm{L}^2$\scriptsize{$\downarrow$} & S-Rate\scriptsize{$\uparrow$} & S-Acc\scriptsize{$\downarrow$} & $\mathrm{L}^2$\scriptsize{$\downarrow$} & S-Rate\scriptsize{$\uparrow$} & S-Acc\scriptsize{$\downarrow$} \\
\hline\hline
\multirow{3}{*}{\rotatebox[origin=c]{90}{\textbf{Linear}}} & Diff-PINN  & \textbf{0.87} &\textbf{100} & \textbf{0.02} & \textbf{0.45} & - & - \\
& VPINN & 0.91 & \textbf{100} & 0.03 & 0.63 & - & - \\
& WPINN & 1.13 & 75.34 & 0.28 & 0.92 & - & - \\
& \tiny \makecell{ WEPINN \\w/o Entropy}   & 0.93 & \textbf{100} & 0.03 & 0.87 & - & - \\
\hline
\multirow{4}{*}{\rotatebox[origin=c]{90}{\textbf{Burgers}}} & Diff-PINN  & 8.62 & 40.75 & 0.79 & 6.23 & 40.00 & \textbf{0.01} \\
& VPINN & 5.34 & 56.13 & 0.43 & 5.89 & 48.90 & 0.05 \\
& WPINN & 9.74 & 92.20 & 0.62 & 7.36 & 60.41 & 0.19 \\
& \tiny \makecell{ WEPINN \\w/o Entropy}   & 3.68 & 82.42 & 0.06 & 3.12 & 99.36 & 0.07 \\
& WEPINN   & \textbf{1.84} & \textbf{98.58} & \textbf{0.05} & \textbf{1.53} & \textbf{99.41} & 0.04 \\
\hline
\multirow{4}{*}{\rotatebox[origin=c]{90}{\textbf{LWR}}} & Diff-PINN  & 6.12 & 41.60 & 1.49 & 3.37 & 32.07 & 0.81 \\
& VPINN & 6.53 & 47.32 & 0.95 & 3.15 & 31.76 & 0.15 \\
& WPINN & 9.72 & 41.45 & 1.48 & 6.36 & 46.23 & 0.78 \\
& \tiny \makecell{ WEPINN \\w/o Entropy}   & 2.54 & 77.22 & 0.19 & 2.50 & 96.45 & 0.28 \\
& WEPINN & \textbf{1.89} & \textbf{93.90} & \textbf{0.15} & \textbf{2.05} & \textbf{96.71} & \textbf{0.23} \\
\hline
\end{tabular}
\caption{Performance comparison under the Dirichlet boundary condition across problems and initial conditions. Reported metrics are the relative $\mathrm{L}^2$ error, the shock detection rate (S-Rate), and the shock position accuracy (S-Acc). Arrows indicate whether lower ($\downarrow$) or higher ($\uparrow$) is better. All numbers are reported in e-2.}
\label{table:dbc}
\end{table}

\begin{table*}[ht]
\centering
\scriptsize
\begin{adjustbox}{center}
\begin{tabular}{|c|c|ccc|ccc|ccc|ccc|}
\hline
\scriptsize{\textbf{Problem}} & \textbf{Method} & \multicolumn{3}{c|}{\textbf{Trig}} & \multicolumn{3}{c|}{\textbf{Fourier}} & \multicolumn{3}{c|}{\textbf{Bell}} & \multicolumn{3}{c|}{\textbf{PWC}} \\
\cline{3-14}
&& $\mathrm{L}^2$ \scriptsize{$\downarrow$}  & S-Rate\scriptsize{$\uparrow$}  & S-Acc\scriptsize{$\downarrow$}   & $\mathrm{L}^2$\scriptsize{$\downarrow$}   & S-Rate\scriptsize{$\uparrow$}  & S-Acc\scriptsize{$\downarrow$}   & $\mathrm{L}^2$\scriptsize{$\downarrow$}   & S-Rate\scriptsize{$\uparrow$}  & S-Acc\scriptsize{$\downarrow$}  & $\mathrm{L}^2$\scriptsize{$\downarrow$}   & S-Rate\scriptsize{$\uparrow$}  & S-Acc\scriptsize{$\downarrow$}  \\
\hline\hline
\multirow{3}{*}{\rotatebox[origin=c]{90}{\textbf{Linear}}}
& Diff-PINN  & 2.05 & - & - & 1.33 & - & - & 1.55 & - & - & 2.99 & 86.09 & 0.30 \\
& VPINN & 2.14 & - & - & 1.30 & - & - & 1.36 & - & - & 1.58 & 85.99 & 0.31 \\
& \tiny \makecell{ WEPINN \\w/o Entropy}  & \textbf{1.89} & - & - & \textbf{0.53} & - & - & \textbf{0.95} & - & - & \textbf{0.96} & \textbf{99.19} & \textbf{0.20} \\
\hline
\multirow{4}{*}{\rotatebox[origin=c]{90}{\textbf{Burgers}}}
& Diff-PINN  & 8.22 & 1.04 & 0.12 & 5.51 & 1.66 & 0.18 & 6.06 & 9.69 & 0.45 & 6.36 & 12.87 & 0.55 \\
& VPINN & 8.47 & 0.45 & 0.05 & 5.26 & 2.57 & 0.10 & 7.31 & 5.86 & 0.08 & 6.15 & 15.78 & 0.52 \\
& \tiny \makecell{ WEPINN \\w/o Entropy}  & 2.15 & 46.57 & \textbf{0.08} & 0.89 & \textbf{95.88} & 0.25 & 0.91 & \textbf{99.85} & 0.23 & 1.64 & \textbf{97.94} & 0.27 \\
& WEPINN  & \textbf{1.80} & \textbf{57.83} & 0.10 & \textbf{0.72} & 94.25 & \textbf{0.11} & \textbf{0.82} &99.40 & \textbf{0.15} & \textbf{1.28} & 95.16 & \textbf{0.22} \\
\hline
\multirow{4}{*}{\rotatebox[origin=c]{90}{\textbf{LWR}}}
& Diff-PINN  & 9.39 & 0.60 & 0.00 & 6.12 & 0.00 & 0.00 & 9.27 & 3.58 & 0.22 & 11.15 & 16.19 & 0.62 \\
& VPINN & 7.50 & 2.32 & 0.14 & 4.57 & 2.94 & 0.40 & 12.00 & 0.95 & 0.23 & 10.96 & 30.59 & 0.75 \\
& \tiny \makecell{ WEPINN \\w/o Entropy}  & 3.36 & \textbf{76.39} & 0.40 & 1.63 & 83.07 & 0.61 & 2.58 & 91.94 & 0.32 & 4.35 & \textbf{98.30} & 0.32 \\
& WEPINN   & \textbf{2.98} & 71.58 & \textbf{0.13} & \textbf{0.80} & \textbf{88.71} & \textbf{0.15} & \textbf{0.95} & \textbf{99.21} & \textbf{0.18} & \textbf{1.85} & 97.78 & \textbf{0.21} \\
\hline
\end{tabular}
\end{adjustbox}
\caption{Performance comparison under the periodic boundary condition across problems and initial conditions. Reported metrics are the relative $\mathrm{L}^2$ error, the shock detection rate (S-Rate), and the shock position accuracy (S-Acc). Arrows indicate whether lower ($\downarrow$) or higher ($\uparrow$) is better. All numbers are reported in e-2.}
\label{table:pbc}
\end{table*}

Tables \ref{table:dbc} and \ref{table:pbc} summarize the results under Dirichlet and periodic boundary conditions, respectively. For the linear advection equation, no shocks form from the Sigmoid initial conditions, so we omit the shock-resolution metrics. In addition, since the entropy condition plays no role in solution selection for linear problems \citep{dafermos1983hyperbolic}, we report only WEPINN without the entropy loss (WEPINN w/o Entropy) for the linear advection equation.

Under the Dirichlet boundary condition, WEPINN achieves substantially lower relative $\mathrm{L}^2$ errors than the baselines on Burgers' equation and the LWR model, together with near-perfect S-Rate and the smallest S-Acc. These gains are consistent across both Sigmoid and Riemann initial conditions.
Both Diff-PINN and VPINN degrade sharply on Burgers' equation and the LWR model, where shocks form. VPINN improves slightly on Diff-PINN in shock localization owing to the use of the weak formulation, but its S-Rate remains too low for S-Acc to be meaningful. WPINN, despite using the weak formulation with the entropy condition, also underperforms on these problems.
By contrast, on the linear advection equation all methods achieve comparable accuracy. In particular, Diff-PINN achieves the lowest relative $\mathrm{L}^2$ error on the smooth Sigmoid initial conditions. This suggests that the failure of Diff-PINN on nonlinear problems stems from the ill-posed pointwise residual near discontinuities rather than from any limitation of the neural network's approximation capacity. These findings show that the proposed WEPINN is the most accurate, particularly in shock resolution, among all methods for nonlinear problems involving shocks.

\begin{figure}[htbp]
    \centering
    \includegraphics[width=1\linewidth]{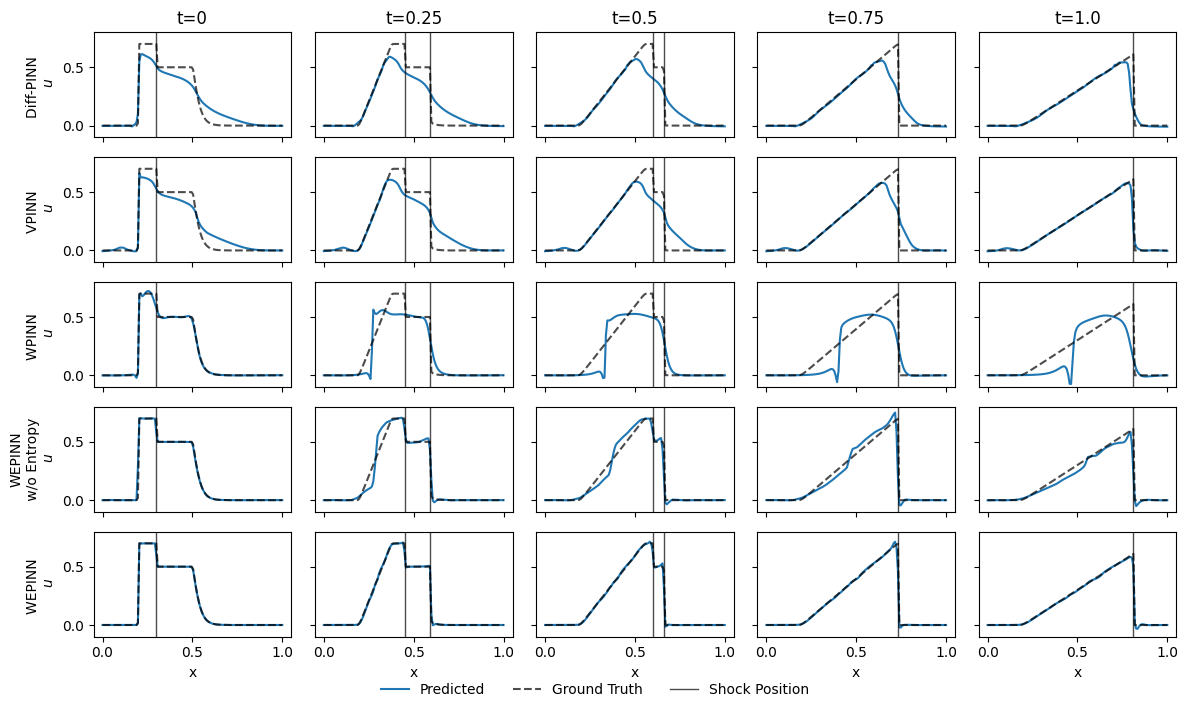}
    \caption{Qualitative comparison on a Burgers' equation example whose solution contains rich phenomena including formation of shocks and rarefactions as well as their interactions. Columns show the solution at $t=0,0.25,0.5,0.75,1.0$. Rows correspond to Diff-PINN, VPINN, WPINN, WEPINN without the entropy loss, and the full WEPINN, together with the ground truth solution in dashed line. Solid vertical lines mark the ground-truth shock locations.}
    \label{fig:combine_compire}
\end{figure}

Under the periodic boundary condition, WEPINN continues to outperform the baselines on all problems. In Table~\ref{table:pbc}, the best result in every cell, across the relative $\mathrm{L}^2$ error, S-Rate, and S-Acc metrics and across all problems and initial conditions, is obtained by either WEPINN or WEPINN w/o Entropy.
The two variants alternate in achieving the best S-Rate and S-Acc, depending on the initial condition class. This suggests that the weak formulation already captures most of the shock structure accurately, while the entropy loss mainly serves to select the physically admissible solution.
Among the baselines, Diff-PINN and VPINN again degrade sharply once shocks form, with S-Rate values below 20 across most initial conditions on Burgers' equation and the LWR model. These findings indicate that the proposed WEPINN is also the most accurate method under the periodic boundary condition. 

To further illustrate the qualitative differences among the methods, we present a representative Burgers' equation example in Figure~\ref{fig:combine_compire}. The example is chosen so that, within a single solution, one observes the formation and propagation of shocks and rarefactions, shock–shock merging, and shock–rarefaction interaction. The columns show the solution at different time snapshots from $t=0$ to $t=1$, while the rows correspond to different methods. The ground-truth solution is shown as a dashed line for comparison, and the ground-truth shock locations are indicated by solid vertical lines.

We observe from Figure~\ref{fig:combine_compire} that Diff-PINN smooths out the solution throughout the temporal evolution. A rarefaction from the left-hand initial discontinuity is captured with noticeable error in magnitude, but the shock formation and propagation, shock-shock merging, and shock-rarefaction interaction are all mispredicted. VPINN offers only a marginal improvement and additionally introduces small spurious oscillations near the boundary. The WPINN prediction shows only a single traveling wave from $t=0.25$ onward, missing most of the wave formation and interaction phenomena. Using our weak-formulation loss alone (WEPINN w/o Entropy), the model accurately captures shock formation and propagation as well as shock–shock merging, but fails to capture the rarefaction. Only the full WEPINN simultaneously captures all the key phenomena in the nonlinear dynamics, including the formation and propagation of shocks and rarefactions, shock–shock merging, and shock–rarefaction interaction, in close agreement with the ground truth.
These observations are consistent with the quantitative results in Tables~\ref{table:dbc} and~\ref{table:pbc}, and underscore the role of the entropy condition: without it, the model cannot distinguish a physical rarefaction from a spurious shock satisfying the weak formulation.
This point is further examined in the ablation study (Section~\ref{sec:role_entropy_loss}), where a simpler example isolates the effect of the entropy loss.

\subsection{One-dimensional system of conservation laws}
\label{sec:vector-1d}
In this subsection, we consider the one-dimensional compressible Euler equations,
\begin{equation} \label{eq:euler-1d}
\partial_t \mathbf{U} + \partial_x \mathbf{F}(\mathbf{U}) = 0,
\end{equation}
where \(\mathbf{U} = (\rho, \rho u, E)^\top\) denotes the vector of conserved quantities, corresponding to mass, momentum, and energy conservation, respectively.
The physical unknowns are the mass density $\rho$, the velocity $u$, and the pressure $p$; the energy density $E$ is determined by the ideal-gas equation of state:
\[
E = \frac{p}{\gamma - 1} + \tfrac{1}{2}\rho u^2
\]
with an adiabatic index \(\gamma > 1\).
The flux function is
\[
\mathbf{F}(\mathbf{U}) =
\begin{pmatrix}
\rho u \\
\rho u^2 + p \\
u(E + p)
\end{pmatrix}.
\]

The Euler system is a hyperbolic system of conservation laws and admits rich wave phenomena, including shocks, contact discontinuities, and rarefactions. As with the scalar case, an entropy condition is required to select the physically admissible solution. For the Euler system, the natural entropy-entropy flux pair is derived from the thermodynamic entropy. Following the classical construction \citep{dafermos1983hyperbolic}, we take
\[
\eta(\vecU) = -\rho s, \quad q(\vecU) = -\rho u s \quad \text{where} \quad s=\ln\!\left(\frac{p}{\rho^{\gamma}}\right).
\]
The pair $(\eta,q)$ satisfies the compatibility condition \eqref{eq:entropy_pair}, and $\eta$ is convex on the physically admissible region $p,\rho>0$.

We evaluate the proposed method on the Sod shock tube problem, a standard benchmark for the compressible Euler equations. From a piecewise-constant initial condition separating two constant states, the solution develops three characteristic wave structures: a rarefaction, a contact discontinuity, and a shock, each propagating at a different speed. Dirichlet boundary conditions are imposed on the domain $\Omega=[0,1]$.

\begin{table}[htbp]
\centering
\begin{tabular}{|c|ccc|ccc|ccc|}
\hline
 Method & \multicolumn{3}{c|}{density $\rho$} & \multicolumn{3}{c|}{speed $u$} & \multicolumn{3}{c|}{pressure $p$} \\ \cline{2-10} 
 & $\mathrm{L}^2$\scriptsize{$\downarrow$} & S-Rate\scriptsize{$\uparrow$} & S-Acc\scriptsize{$\downarrow$}  & $\mathrm{L}^2$\scriptsize{$\downarrow$} & S-Rate\scriptsize{$\uparrow$} & S-Acc\scriptsize{$\downarrow$}  &  $\mathrm{L}^2$\scriptsize{$\downarrow$} & S-Rate\scriptsize{$\uparrow$} & S-Acc\scriptsize{$\downarrow$}  \\ \hline \hline
Diff-PINN & 6.93 & 37 & 2.45 & 32.2 & 0 & N.\,A. & 3.46 & 25 & 0.35 \\ \hline
VPINN & 6.53 & 43 & 1.88 & 33.7 & 5 & 1.95 & 3.15 & 29 & 0.32 \\ \hline
WPINN & 1.38 & - & - & 3.77 & - & - & 1.47 & - & - \\ \hline
\scriptsize \makecell{WEPINN \\w/o Entropy} & 0.90 & 79 & 0.52 & 2.87 & 99 & 0.51 & 1.13 & 100 & 0.36 \\ \hline
WEPINN & \textbf{0.60} & \textbf{99} & \textbf{0.30} & \textbf{2.34} & \textbf{100} & \textbf{0.37} & \textbf{0.61} & \textbf{100} & \textbf{0.29} \\ \hline
\end{tabular}
\caption{Performance comparison on the Sod shock tube problem for the compressible Euler equations. Reported metrics are the relative $\mathrm{L}^2$ error, the shock detection rate (S-Rate), and the shock position accuracy (S-Acc). Arrows indicate whether lower ($\downarrow$) or higher ($\uparrow$) is better. All numbers are reported in e-2.}
\label{tab:euler}
\end{table}

\begin{figure}[htbp]
    \centering
    \includegraphics[width=1\linewidth]{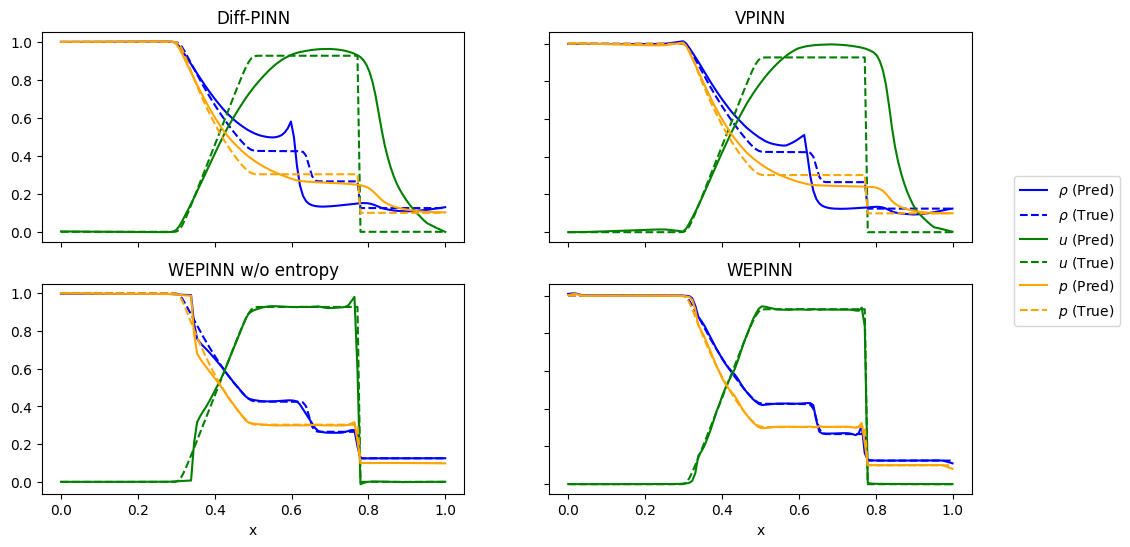}
    \caption{Comparison of predicted solutions (density $\rho$, speed $u$, and pressure $p$ in different colors) at $t=0.2$ for the Sod shock tube problem for the compressible Euler equations. Solid lines indicate predictions and dashed lines indicate the ground truth. Each panel shows a different method. }
    \label{fig:euler}
\end{figure}

Table~\ref{tab:euler} reports the relative $\mathrm{L}^2$ error, S-Rate, and S-Acc of each method on the density $\rho$, velocity $u$, and pressure $p$ predicted by the models, and Figure~\ref{fig:euler} visualizes the predicted solutions at time $t=0.2$. For WPINN, we quote the numbers reported in \citet{chaumet2022improving}; as discussed in Section~\ref{sec:setup}, no implementation is available to us for the Euler equations, so S-Rate, S-Acc, and the solution visualization are absent.

We observe from Table~\ref{tab:euler} that WEPINN achieves the lowest $\mathrm{L}^2$ error on all three variables, together with the highest S-Rate and lowest S-Acc. Figure~\ref{fig:euler} shows that WEPINN accurately reproduces all three waves: the rarefaction on the left, the contact discontinuity in the middle, and the shock on the right. Removing the entropy loss (WEPINN w/o Entropy) preserves the shock and the contact discontinuity but yields a less accurate rarefaction, consistent with the behavior observed on the scalar problems in Section~\ref{sec:scalar_1d}.
Diff-PINN and VPINN perform comparably, with VPINN offering only a marginal improvement, and both methods yield smeared profiles across the spatial domain. These results confirm the effectiveness of WEPINN on this benchmark problem of compressible Euler equations, and that entropy enforcement remains essential in the presence of more complex wave structures.

\subsection{Two-dimensional scalar conservation laws}
\label{sec:scalar-2d}
In this subsection, we consider a scalar Burgers' equation in two space dimensions, obtained by taking $f(u)=\frac12u^2$ as the flux in each spatial direction:
\begin{equation}\label{eq:burgers2d}
\partial_t u + \grad \cdot \vecF(u) = 0 \quad \text{with} \quad \vecF(u) = \left(\frac{u^2}{2},\ \frac{u^2}{2}\right),
\end{equation}
or equivalently
\begin{equation}
\partial_t u + \partial_x \left(\frac{u^2}{2}\right) + \partial_y \left(\frac{u^2}{2}\right) = 0.
\end{equation}
This is a special case of the two-dimensional scalar conservation laws studied in \citet{zhang1989two}, and a natural two-dimensional extension of the one-dimensional Burgers' equation in Section~\ref{sec:scalar_1d}.

We consider the following two initial conditions:
\begin{itemize}
    \item a \emph{disk} profile
    \begin{equation*}
u_0(x,y)=
\begin{cases}
0.5, & (x-0.5)^2+(y-0.5)^2 < 0.04;\\
0, & \text{otherwise}
\end{cases}
\end{equation*}
which is piecewise constant with a circular discontinuity of radius $0.2$ centered at $(0.5,0.5)$;
\item a smooth \emph{trigonometric} profile \begin{equation*}
u_0(x,y)=\frac{1}{4}\sin(2\pi x)\sin(2\pi y).
\end{equation*}
\end{itemize}
The two initial conditions test complementary regimes: the disk one contains a circular initial discontinuity whose propagation needs to be tracked, while the trigonometric one is smooth and tests the model's ability to capture shock formation from smooth data.
We impose periodic boundary conditions in both spatial directions on the domain $\Omega=[0,1]^2$. WPINN is not included in the comparison because no implementation for two-dimensional problems is available, as discussed in Section~\ref{sec:setup}.

\begin{table}[htbp]
\centering
\begin{tabular}{|c|c|c|}
\hline
Method & Disk & Trigonometric \\ \hline\hline
Diff-PINN & 0.28 & 0.17 \\ \hline
VPINN & 0.27 & 0.21 \\ \hline
\scriptsize \makecell{ WEPINN \\w/o Entropy} & 0.96 & 0.57 \\ \hline
WEPINN & \textbf{0.12} & \textbf{0.14} \\ \hline
\end{tabular}
\caption{Relative $\mathrm{L}^2$ errors on the two-dimensional scalar Burgers' equation under periodic boundary conditions.}
\label{table:2DBurger}
\end{table}

\begin{figure}[htbp]
    \centering
    \includegraphics[width=0.45\linewidth]{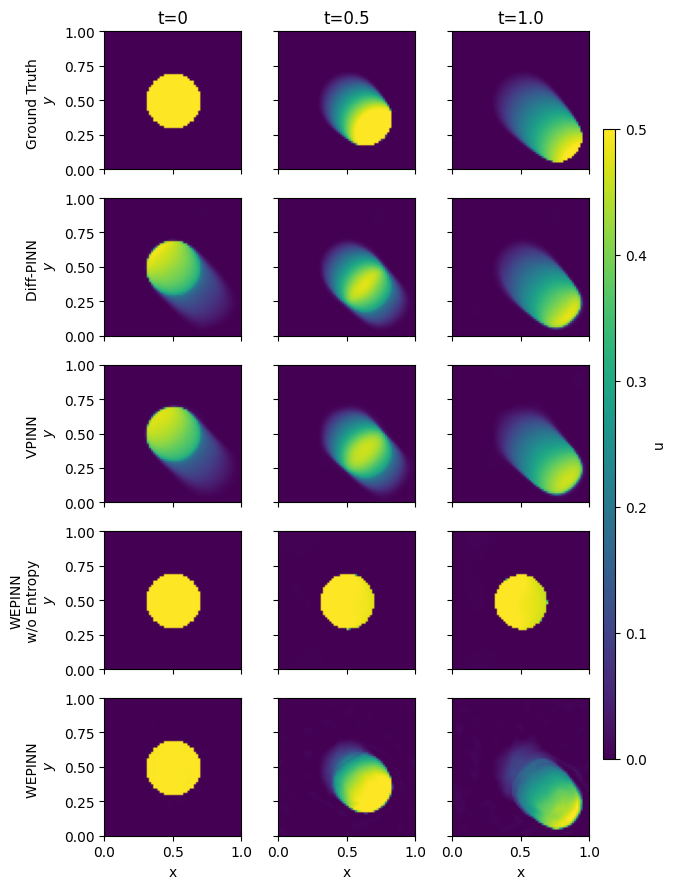}
    \includegraphics[width=0.45\linewidth]{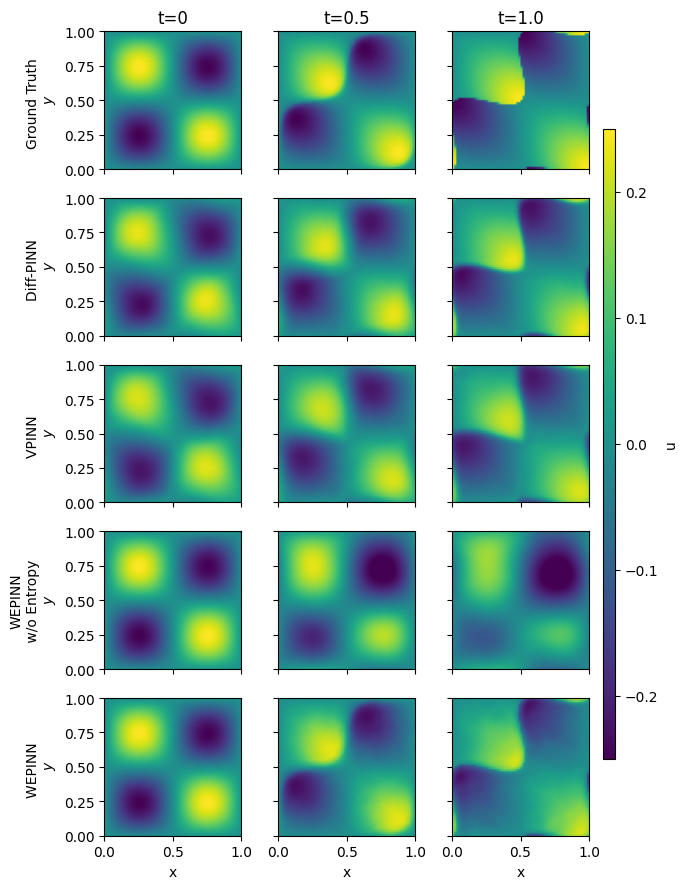}
    \caption{Predicted solutions of the two-dimensional scalar Burgers' equation under periodic boundary conditions. (Left) Disk initial condition; (Right) Trigonometric initial condition. Within each panel, rows correspond to the ground truth and different methods. Columns show the solution at $t=0$, $t=0.5$, and $t=1$.}
    \label{fig:Burger2D}
\end{figure}

Table~\ref{table:2DBurger} reports the relative $\mathrm{L}^2$ error of each method on the two initial conditions, and Figure~\ref{fig:Burger2D} visualizes the predicted solutions at $t=0$, $t=0.5$, and $t=1$. Because shocks in two dimensions are curves rather than isolated points, the S-Rate and S-Acc metrics from Section~\ref{sec:setup} do not directly extend to this setting; we therefore report only the $\mathrm{L}^2$ errors.

Across both initial conditions, WEPINN achieves the lowest $\mathrm{L}^2$ error. Moreover, Figure~\ref{fig:Burger2D} shows that WEPINN accurately tracks the propagation of the circular discontinuity in the disk case and captures shock formation from the smooth trigonometric profile. By contrast, Diff-PINN and VPINN smooth out the solution on both initial conditions, particularly at the intermediate time steps, and in some cases distort the initial condition itself.

The behavior of WEPINN without the entropy loss (WEPINN w/o Entropy) reveals a notable difference from the experiments on the one-dimensional scalar problems (Section~\ref{sec:scalar_1d}), where WEPINN w/o Entropy is already competitive with the full WEPINN. In the two-dimensional setting, however, it produces the largest $\mathrm{L}^2$ error among all methods and fails to reproduce shock formation and propagation. This contrast suggests that the entropy loss becomes increasingly important in higher-dimensional problems.

\section{Ablation Studies}\label{sec:ablation}

\subsection{Role of the entropy loss}\label{sec:role_entropy_loss}

\begin{figure}[htbp]
    \centering
    \includegraphics[width=0.8\linewidth]{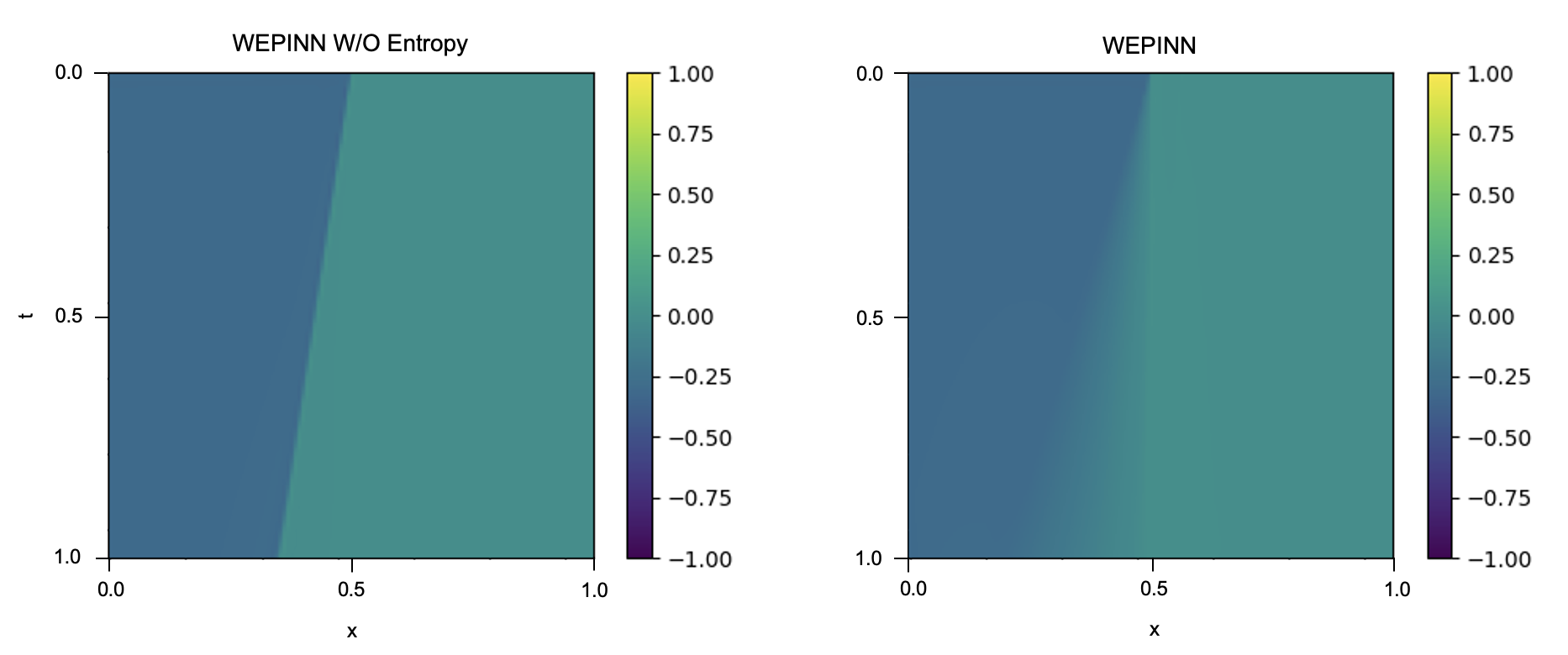}
    \caption{Role of the entropy loss in WEPINN, illustrated on the Riemann problem for Burgers' equation \eqref{eq:ablation_burgers} with $x_0=0.5$, $u_L= -0.2$, and $u_R=0$, under the Dirichlet boundary condition. (Left) WEPINN without the entropy loss; (Right) Full WEPINN. In each panel, the horizontal axis is the spatial coordinate $x$ and the vertical axis is time $t$.}
    \label{fig:Entropy_role}
\end{figure}

The entropy condition \eqref{eq:entropy} is essential for eliminating unphysical weak solutions to the hyperbolic conservation law \eqref{eq:HCL} that satisfy the weak formulation \eqref{eq:weak} but violate physical admissibility.
A classical illustrative example is the Riemann problem for the one-dimensional Burgers' equation
\begin{equation}\label{eq:ablation_burgers}
\partial_tu+\partial_x\left(\frac{u^2}2\right)=0, \qquad u(0,x)=\begin{cases}
    u_L, \quad x<x_0 \\ u_R, \quad x>x_0
\end{cases}
\end{equation}
with $u_L<u_R$. For this problem, both a shock and a rarefaction satisfy the weak formulation, but only the rarefaction is admissible under the entropy condition.
In WEPINN, the entropy loss enforces the entropy condition during training, so that the model converges to the weak entropy solution.
Figure~\ref{fig:Entropy_role} illustrates this on the above Riemann problem with $x_0=0.5$, $u_L= -0.2$, and $u_R=0$, under the Dirichlet boundary condition. Without the entropy loss, WEPINN converges to a propagating shock, while adding the entropy loss recovers the correct rarefaction. This sheds light on why the full WEPINN outperforms WEPINN without the entropy loss in most cases across the numerical experiments in Section~\ref{sec:results}, particularly on initial conditions where rarefactions develop.

\subsection{Test function selection}\label{sec:test-function-abl}

The choice of test functions used in the weak-entropy loss is critical to the performance of WEPINN.
This includes both the function family (e.g., trigonometric functions or orthogonal polynomials such as Chebyshev polynomials) and the maximum degree of the functions, which controls their expressive power.

To determine the optimal maximum degree (i.e., maximum frequency) of the trigonometric test functions chosen in Section~\ref{eq:method_wepinn}, we conduct an ablation study by varying 
the degree from 4 to 32. 
The results are presented in Table~\ref{table:degree}. 
We observe that degrees 16 and 32 yield the best performance. 
Moreover, degree 16 performs better on high-frequency data, while degree 32, which is used in the WEPINN implementation for the numerical experiments in Section~\ref{sec:results}, yields better results in scenarios with sharp discontinuities. 

\begin{table*}[ht]
\centering
\scriptsize
\begin{adjustbox}{center}
\begin{tabular}{|c|c|ccc|ccc|ccc|ccc|}
\hline
\scriptsize{\textbf{Problem}} & \textbf{Degree} & \multicolumn{3}{c|}{\textbf{Trig}} & \multicolumn{3}{c|}{\textbf{Fourier}} & \multicolumn{3}{c|}{\textbf{Bell}} & \multicolumn{3}{c|}{\textbf{PWC}} \\
\cline{3-14}
&& $\mathrm{L}^2$ \scriptsize{$\downarrow$}  & S-Rate\scriptsize{$\uparrow$}  & S-Acc\scriptsize{$\downarrow$}   & $\mathrm{L}^2$\scriptsize{$\downarrow$}   & S-Rate\scriptsize{$\uparrow$}  & S-Acc\scriptsize{$\downarrow$}   & $\mathrm{L}^2$\scriptsize{$\downarrow$}   & S-Rate\scriptsize{$\uparrow$}  & S-Acc\scriptsize{$\downarrow$}  & $\mathrm{L}^2$\scriptsize{$\downarrow$}   & S-Rate\scriptsize{$\uparrow$}  & S-Acc\scriptsize{$\downarrow$}  \\
\hline\hline
\multirow{4}{*}{\rotatebox[origin=c]{90}{\textbf{Burgers}}}
& 4 & 9.27 & 16.35 & 0.23 & 4.52 & 24.11 & 0.97 & 3.32 & 73.88 & 1.03 & 3.80 & 82.98 & 0.85 \\
& 8 & 2.91 & 40.22 & 0.16 & 1.85 & 52.76 & 0.33 & 2.43 & 96.43 & 0.42 & 2.06 & 87.51 & 0.24 \\
& 16 & \textbf{1.43} & 54.69 & 0.16 & 0.96 & 65.84 & 0.14 & 1.27 & 98.94 & 0.18 & \textbf{1.21} & 95.06 & 0.23\\
& 32 & 1.80 & \textbf{57.83} & \textbf{0.10} & \textbf{0.72} & \textbf{94.25} & \textbf{0.11} & \textbf{0.82} & \textbf{99.40} & \textbf{0.15} & 1.28 & \textbf{95.16} & \textbf{0.22} \\
\hline
\multirow{4}{*}{\rotatebox[origin=c]{90}{\textbf{LWR}}}
& 4 & 10.14 & 22.26 & 0.20 & 3.15 & 47.77 & 0.60 & 4.90 & 59.72 & 0.63 & 5.09 & 90.04 & 0.39 \\
& 8 & 4.38 & 56.08 & 0.68 & 1.53 & 76.42 & 0.39 & 2.27 & 79.71 & 0.41 & 3.05 & 91.70 & 0.32 \\
& 16 & 3.06 & 63.93 & 0.31 & 0.97 & 83.69 & 0.24 & 1.94 & 91.43 & 0.19 & 2.24 & 96.73 & 0.24\\
& 32   & \textbf{2.98} & \textbf{71.58} & \textbf{0.13} & \textbf{0.80} & \textbf{88.71} & \textbf{0.15} & \textbf{0.95} & \textbf{99.21} & \textbf{0.18} & \textbf{1.85} & \textbf{97.78} & \textbf{0.21} \\
\hline
\end{tabular}
\end{adjustbox}
\caption{Performance comparison between different maximum degrees of trigonometric test functions in WEPINN across different problems and initial conditions under the periodic boundary condition. Reported metrics are the relative $\mathrm{L}^2$ error, the shock detection rate (S-Rate), and the shock position accuracy (S-Acc). Arrows indicate whether lower ($\downarrow$) or higher ($\uparrow$) is better. All numbers are reported in e-2.}
\label{table:degree}
\end{table*}

\begin{table}[htbp]
\centering
\footnotesize
\setlength{\tabcolsep}{4pt}
\renewcommand{\arraystretch}{1.2}
\begin{tabular}{|c|c|ccc|ccc|}
\hline
\scriptsize{\textbf{Problem}} & \textbf{Test function} & \multicolumn{3}{c|}{\textbf{Riemann}} & \multicolumn{3}{c|}{\textbf{Sigmoid}} \\
\cline{3-8}
&& $\mathrm{L}^2$\scriptsize{$\downarrow$} & S-Rate\scriptsize{$\uparrow$} & S-Acc\scriptsize{$\downarrow$} & $\mathrm{L}^2$\scriptsize{$\downarrow$} & S-Rate\scriptsize{$\uparrow$} & S-Acc\scriptsize{$\downarrow$} \\
\hline\hline
\multirow{2}{*}{\textbf{Linear}} 
& Cheby & 4.58 & 88.48 & 1.11 & 5.29 & - & - \\
& Trig  & \textbf{0.93} & \textbf{100} & \textbf{0.03} & \textbf{0.87} & - & - \\
\hline
\multirow{2}{*}{\textbf{Burgers}} 
& Cheby  & 44.11 & 77.47 & 1.35 & 46.29 & 73.86 & 1.43 \\
& Trig   & \textbf{1.84} & \textbf{98.58} & \textbf{0.05} & \textbf{1.53} & \textbf{99.41} & \textbf{0.04} \\
\hline
\multirow{2}{*}{\textbf{LWR}} 
& Cheby  & 43.21 & 83.48 & 1.07 & 43.80 & 79.19 & 1.43\\
& Trig & \textbf{1.89} & \textbf{93.90} & \textbf{0.15} & \textbf{2.05} & \textbf{96.71} & \textbf{0.23} \\
\hline
\end{tabular}
\caption{Performance comparison between trigonometric functions and Chebeshev polynomials as test functions in WEPINN across different problems and initial conditions under the periodic boundary condition. Reported metrics are the relative $\mathrm{L}^2$ error, the shock detection rate (S-Rate), and the shock position accuracy (S-Acc). Arrows indicate whether lower ($\downarrow$) or higher ($\uparrow$) is better. All numbers are reported in e-2.}
\label{table:cheby}
\end{table}

To assess the effectiveness of trigonometric test functions, we implemented WEPINN using Chebyshev polynomials of the first kind to construct test functions. For a fair comparison, the maximum degree of the Chebyshev polynomials is set to 32, matching the trigonometric case. Due to the non-periodic nature of Chebyshev polynomials, this comparison is performed under the Dirichlet boundary condition. The results are presented in Table~\ref{table:cheby}. 

We observe from Table~\ref{table:cheby} that using Chebyshev polynomials for test functions results in worse performance than trigonometric functions across all metrics.
In particular, the S-Acc degrades more significantly than S-Rate, suggesting that the model with Chebyshev polynomials struggles to track shocks accurately. 
This may due to the increasing challenge to evaluate integrals in the weak formulation and the entropy condition with Chebeshev polynomials that exhibit rapid oscillations near the domain boundary.

\section{Conclusion}\label{sec:conclusion}
In this work, we introduce a novel physics-informed learning framework that leverages weak entropy conditions to effectively solve hyperbolic conservation laws, where solutions may develop discontinuities. By rigorously adhering to weak formulation and entropy conditions, and employing the discrete fast Fourier transform (DFFT) for efficient numerical quadrature, our method avoids common pitfalls such as unphysical dissipation and oscillations observed in existing approaches. Our extensive numerical experiments demonstrate that the proposed method consistently outperforms existing approaches across a variety of scalar conservation laws and systems of conservation laws in one and two dimensions. It notably excels in accurately resolving complex nonlinear wave phenomena, including formation and propagation of shocks, shock-shock merging, and shock-rarefaction interaction, under diverse initial and boundary conditions. The enforcement of entropy conditions significantly enhances the reliability of our method, particularly by ensuring physical admissibility of the solution. Future work will explore extending this methodology to design neural operators for solving hyperbolic conservation laws and related PDEs involving discontinuities and moving interfaces.

\section*{Acknowledgments}
This work was supported by the National Science Foundation (NSF) under Grant No. CPS-2038984.

\appendix

\section{Evaluation Metrics}

To quantitatively evaluate the performance of different models in solving hyperbolic conservation laws, we adopt several metrics that assess both global solution accuracy and shock-capturing capability. A reference solution is generated using the Lax--Friedrichs numerical scheme on an ultra-fine grid to serves as the ground truth for all subsequent comparisons.

We compute the \textbf{relative $\mathrm{L}^2$ error} between the predicted solution \(\hat{u}\) and the ground truth \(u\) over the entire spatiotemporal domain: \(\text{Relative } \mathrm{L}^2 = \| \hat{u} - u \|_2/\| u \|_2\).This metric captures the overall discrepancy between the predicted and true solutions.

To evaluate the model's ability to identify shocks, we define the \textbf{shock detection rate (S-Rate)} as the proportion of ground truth shock locations that are successfully identified by the model within a specified tolerance band. A shock is considered detected, if a high-gradient region in the predicted solution is located within a small neighborhood (10\% of spatial domain) of the ground truth shock position.

We also compute the \textbf{shock position accuracy (S-Acc)}, defined as the average distance between detected shock locations and their corresponding ground truth positions. This metric captures the spatial precision of shock prediction and penalizes shifts or smearing of the shock front.

The detailed procedure for calculate both newly proposed metric are as shown in Algorithm~\ref{alg:metric}.

\begin{algorithm}[H]
\caption{Computing Shock Detection Rate and Shock Accuracy}
\begin{algorithmic}[1]
\Function{ShockError}{$u$, $\hat{u}$, $dx$, $t_{\max}$, $h$}
    \State Compute gradient: $g_u \gets |\partial_x u|$, $g_{\hat{u}} \gets |\partial_x \hat{u}|$
    \State Compute total variation $TV(t)$ over time
    \State Compute average total variation: $TV_{\text{mean}} \gets \frac{1}{t_{\max}} \sum_t TV(t)$
    \State Initialize counters: $D \gets 0$, $M \gets 0$, $N \gets 0$
    \For{$t = 0$ to $t_{\max}-1$}
        \State Detect shocks: $S_u \gets$ peaks$(g_u[t], \text{height} = h \cdot TV_{\text{mean}})$
        \State Detect shocks: $S_{\hat{u}} \gets$ peaks$(g_{\hat{u}}[t], \text{height} = h \cdot TV_{\text{mean}})$
        \State Compute distance matrix $d_{ij} = \min(|S_u[i] - S_{\hat{u}}[j]|, 128 - |S_u[i] - S_{\hat{u}}[j]|)$
        \For{each $i$ in $S_u$}
            \State $j^* \gets \arg\min_j d_{ij}$
            \If{$d_{ij^*} < d_{\text{max}}$}
                \State $D \gets D + d_{ij^*}$ \Comment{Accumulate distance}
                \State $M \gets M + 1$ \Comment{Count matched shocks}
            \EndIf
            \State $N \gets N + 1$ \Comment{Count total true shocks}
        \EndFor
    \EndFor
    \State Detection rate: $\text{S-Rate} \gets M / N$
    \State Accuracy (distance): $\text{S-Acc} \gets (D / M) \cdot dx$
    \State \Return $\text{S-Rate}, \text{S-Acc}$
\EndFunction
\end{algorithmic}
\label{alg:metric}
\end{algorithm}

\label{sec:entropy}

\section{Initial and Boundary conditions}
\label{sec:icbc}
\subsection{Initial Conditions}

To ensure robustness to diverse solution behaviors, we consider multiple classes of initial conditions reflecting varying smoothness and discontinuity.

For Dirichlet boundary conditions, two types of initial conditions are considered:
\begin{itemize}
    \item \textbf{Sigmoid-type profiles:} Smooth transitions between prescribed left and right boundary states, with adjustable steepness and transition location. These profiles test the model’s ability to resolve sharp but continuous gradients while respecting fixed boundary values.
    \item \textbf{Riemann-type profiles:} Piecewise-constant initial states composed of two or three constant regions. This setting enables systematic evaluation of shock, rarefaction, and their interactions under Dirichlet boundary conditions, including shock–shock and shock–rarefaction configurations.
\end{itemize}

For periodic boundary conditions, we consider four types of initial conditions are used:
\begin{itemize}
    \item \textbf{Fourier series:} Smooth functions constructed from sine and cosine modes with randomly sampled coefficients up to mode 4.
    \item \textbf{Trigonometric functions:} High-frequency sinusoidal waves (up to frequency 4) with random phase shifts, testing sensitivity to oscillatory dynamics.
    \item \textbf{Bell-shaped functions:} Smooth, single-peaked profiles with random centers and widths, designed to assess shock formation from localized features.
    \item \textbf{Piecewise constant functions:} Discontinuous profiles with randomly sampled segment values and breakpoints, modeling shock-like initial states.
\end{itemize}

Periodic boundary conditions enforce continuity across the domain boundaries,
\[
u(t,0) = u(t,1),
\]
modeling a spatially repeating (closed-loop) system. This setting eliminates artificial boundary effects and enables clean analysis of wave propagation and interaction, making it particularly suitable for Fourier-based representations and operator learning methods.

\section{Implementation Detail}
\subsection{Neural Network Architecture}

To better capture complex shock interactions, we employ a neural network architecture with residual (skip) connections and ReLU activation. The network uses a hidden dimension of 128 and consists of 9 residual layers. The detailed architecture is presented in Algorithms~\ref{alg:res} and~\ref{alg:respinn}. For baseline comparisons, we replace the ReLU activation with Tanh, as ReLU can negatively affect the automatic differentiation required for computing the differential form loss in PINN.
\begin{algorithm}[H]
\caption{Residual Block (ResBlock)}
\begin{algorithmic}[1]
\Function{ResBlock}{$x$}
    \State $residual \gets x$
    \State $x \gets \text{fc}_1(x)$
    \State $x \gets \text{activation}(x)$
    \State $x \gets \text{fc}_2(x)$
    \State \Return $x + residual$
\EndFunction
\end{algorithmic}
\label{alg:res}
\end{algorithm}

\begin{algorithm}[H]
\caption{PINN with Residual Connections (PINNRes)}
\begin{algorithmic}[1]
\Function{PINNRes}{$x$, $t$}
    \State $input \gets \text{concat}(x, t)$
    \State $x \gets \text{fc}_{in}(input)$
    \State $x \gets \text{activation}(x)$
    \For{$i = 1$ to $L$}
        \State $x \gets \text{ResBlock}(x)$
    \EndFor
    \State $output \gets \text{fc}_{out}(x)$
    \State \Return $output$
\EndFunction
\end{algorithmic}
\label{alg:respinn}
\end{algorithm}

\subsection{Training}
The model is trained using the Adam optimizer with an initial learning rate of $1 \times 10^{-3}$. 
A scheduled learning rate decay is applied, reducing the learning rate by a factor of 0.3 at steps 10{,}000 and 20{,}000. 
The total number of training steps is 30{,}000. All reported numbers in the table are average across 30 samples.

\section{Compatibility of different families of test functions with uniform grids}
\label{sec:test-function-explan}
\begin{figure}
    \centering
    \includegraphics[width=0.45\linewidth]{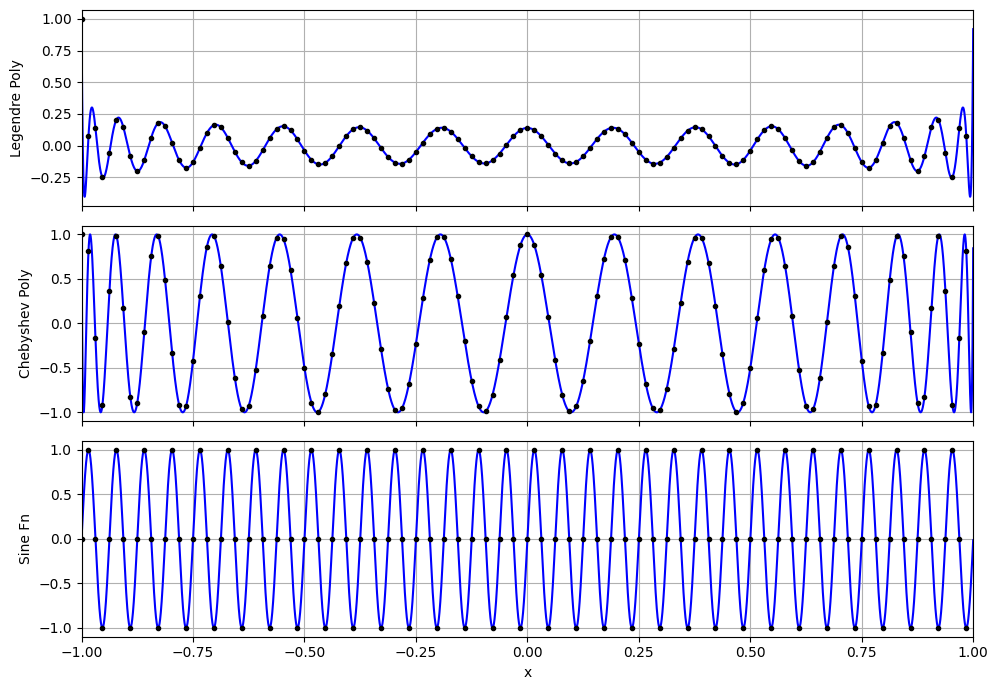}
    \includegraphics[width=0.45\linewidth]{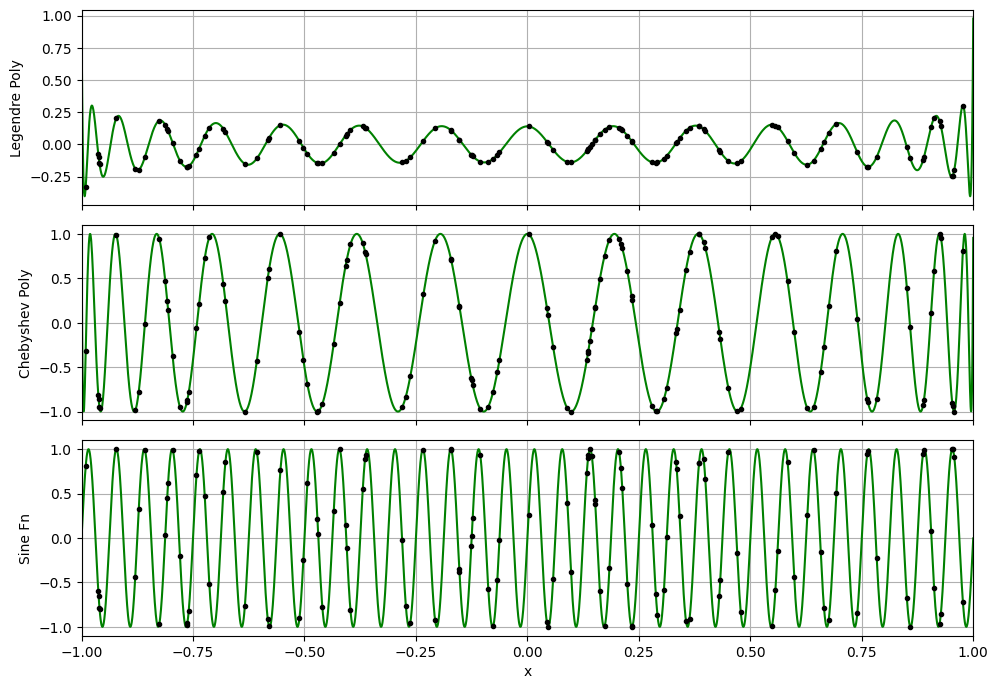}
    \caption{Comparison of evaluations of a Legendre polynomial of degree 32, a Chebyshev polynomial of degree 32, and a trigonometric function of frequency 32 on uniform grid and random collocation points.}
    \label{fig:test_uniform}
\end{figure}

In prior PINN research that incorporates test functions—such as VPINN or WPINN—Legendre or Chebyshev polynomials are commonly used due to their orthogonality properties. While any complete basis that spans the desired test function space is theoretically sufficient, we opt for trigonometric functions instead of polynomials. The key motivation lies in their uniform wavelength and better behavior under discretizations. At higher frequencies, polynomial bases become increasingly difficult to resolve accurately on finite grids, leading to numerical artifacts or instability. In contrast, trigonometric functions maintain consistent oscillatory behavior across frequencies, making them more suitable for capturing fine-scale solution structures. This difference is as illustrated in Figures~\ref{fig:test_uniform} left. To validate this design choice, we additionally implemented our proposed weak and entropy losses using Chebyshev polynomial test functions for comparison in the section~\ref{sec:test-function-abl}.

Similarly, we opt for a uniform grid rather than a random grid that is commonly used in prior PINN works, as a uniform grid has a better capacity to capture high frequency test functions, whereas a random grid could miss peaks in the waves of test functions, as shown in Figure~\ref{fig:test_uniform} right.

\section{Efficiency}
\label{sec:efficiency}

The main computational bottleneck of the proposed method lies in evaluating $\mathcal{I}_{t,\vecx}^{\varphi_n}$ and $\mathcal{J}_{t,\vecx}^{\varphi_n}$, which involve summations over the full space-time grid. Here we take one dimensional scalar case as an example:
\begin{align}
I_{t,x}^\varphi &= \sum_{j=0}^{n_x}\sum_{k=0}^{n_t} \alpha_k \beta_j
\big(u\, \partial_t \varphi + f(u)\, \partial_x \varphi \big)(t_k,x_j), \label{eq:Ixt_app} \\
J_{t,x}^\varphi &= \sum_{j=0}^{n_x}\sum_{k=0}^{n_t} \alpha_k \beta_j
\big(\eta(u)\, \partial_t \varphi + q(u)\, \partial_x \varphi \big)(t_k,x_j). \label{eq:Jxt_app}
\end{align}
Direct evaluation of these summations incurs a computational cost of $O(N^{2d})$, which becomes prohibitive for fine grids.

\paragraph{DFFT acceleration.}
Given choosing trigonometric test functions on a uniform grid, the above summations can be interpreted as discrete Fourier transforms of the weighted solution. This enables efficient computation using the discrete fast Fourier transform (DFFT).

To illustrate, consider the first term in~\eqref{eq:Ixt_app} with a test function
\[
\varphi(t,x) = \sin(2\pi n t)\sin(2\pi m x), \quad n,m \in \mathbb{Z}.
\]
Then,
\begin{align*}
\sum_{j=0}^{n_x}\sum_{k=0}^{n_t} \alpha_k \beta_j u(t_k,x_j) \partial_t \varphi(t_k,x_j)
&=
\sum_{j=0}^{n_x}\sum_{k=0}^{n_t} \alpha_k \beta_j u(t_k,x_j)
\left(2\pi n \cos(2\pi n t_k)\sin(2\pi m x_j)\right).
\end{align*}

Using Euler’s identity, the product of trigonometric functions can be expressed as a linear combination of complex exponentials:
\begin{align*}
\cos(2\pi n t_k)\sin(2\pi m x_j)
=
\frac{1}{2}\Im\left(
e^{2\pi i (n t_k + m x_j)} + e^{2\pi i (-n t_k + m x_j)}
\right).
\end{align*}

Substituting this into the summation yields
\begin{align*}
\sum_{j=0}^{n_x}\sum_{k=0}^{n_t} \alpha_k \beta_j u(t_k,x_j) \partial_t \varphi(t_k,x_j)
=
\pi n \sum_{j=0}^{n_x}\sum_{k=0}^{n_t} \alpha_k \beta_j u(t_k,x_j)
\Im\left(
e^{2\pi i (n t_k + m x_j)} + e^{2\pi i (-n t_k + m x_j)}
\right).
\end{align*}

Each term corresponds to a Fourier coefficient of the weighted field $\alpha_k \beta_j u(t_k,x_j)$. Therefore, all such terms can be computed simultaneously using DFFT.

\paragraph{Computational complexity.}
Using DFFT, the computational cost is reduced to
\[
O(n_t n_x \log(n_t n_x)),
\]
compared to the quadratic cost $O(n_t^2 n_x^2)$ of direct evaluation. In $d$ spatial dimensions, this generalizes to
\[
O(N_tN_x^d \log N_t\log^d N_x),
\]
which is significantly more efficient than the $O(N_t^2N_x^{2d})$ complexity of standard numerical integration.

\paragraph{Practical efficiency.}
This acceleration is particularly important when using high-frequency test functions, which are necessary to capture fine-scale structures in the solution. Without DFFT, evaluating the weak and entropy residuals becomes both time- and memory-intensive. 

\begin{figure}
    \centering
    \includegraphics[width=0.8\linewidth]{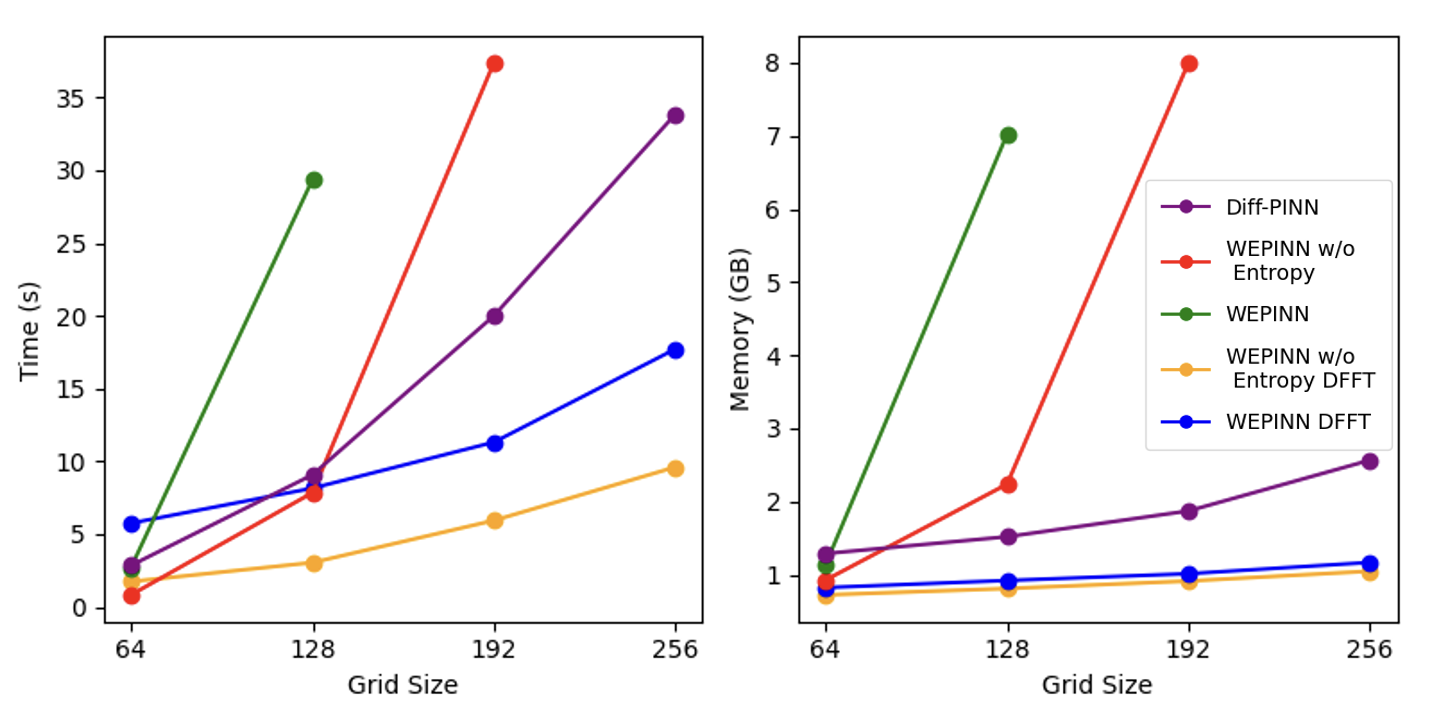}
    \caption{Time and memory efficiency of the proposed method compared against the Diff-PINN and our proposed method with and without DFFT. Left: Runtime (in seconds) for 500 training epochs across varying grid sizes on single Nvidia RTX4090. Right: Total memory consumption by the model during training. The DFFT-based implementation demonstrates superior scalability and resource efficiency.}
    \label{fig:eff}
\end{figure}

We validate the efficiency of the proposed implementation through an ablation study on runtime and memory consumption across different grid resolutions, as shown in Figure~\ref{fig:eff}. Due to GPU memory limitations, direct numerical integration could not be evaluated at larger grid sizes.

\bibliographystyle{plainnat} 
\bibliography{references}

\end{document}